\documentclass[11pt]{article}
\usepackage{latexsym}
\usepackage{amssymb}
\usepackage{graphicx}
\usepackage{cite}
\usepackage{amsmath}
\usepackage{mathrsfs}
\usepackage{color}
\usepackage{subfigure}

\newtheorem{theorem}{Theorem}[part]
\newtheorem{definition}{Definition}[part]

\newtheorem{lemma}{Lemma}[part]
\newtheorem{corollary}{Corollary}[part]
\newtheorem{remark}{Remark}[part]
\newtheorem{example}{Example}[part]

\def \ep{\hbox{ }\hfill$\Box$}

\begin{document}

\title{Optimality Conditions for Sparse Bilinear Least Squares Problems\thanks{This work was partially supported by the National Natural Science Foundation of China under Grants 12371309, 12471295 and 12426306,   Guangdong Basic and Applied Basic Research Foundation (2024A1515011566), and Hetao Shenzhen-Hong Kong Science and Technology  Innovation Cooperation Zone Project (HZQSWS-KCCYB-2024016).} }

\author{Zixin Deng\thanks{School of Mathematics, Tianjin University, Tianjin 300072, China.  Email: dengzx@tju.edu.cn.}
\and
Zheng-Hai Huang\thanks{\emph{Corresponding author.} School of Mathematics, Tianjin University, Tianjin 300072, China. Email: huangzhenghai@tju.edu.cn.}
\and
Yun-Bin Zhao\thanks{Shenzhen International Center for Industrial and Applied Mathematics, SRIBD, The Chinese University of Hong Kong, Shenzhen 518172, China. Email: yunbinzhao@cuhk.edu.cn.}
}

\date{}

\maketitle

\begin{abstract}
The first-order optimality conditions of sparse bilinear least squares problems are studied. The so-called T-type and N-type stationary points for this problem are characterized in terms of tangent cone and normal cone in Bouligand and Clarke senses, and another stationarity concept called the coordinate-wise minima is introduced and discussed. Moreover, the L-like stationary point for this problem is introduced and analyzed through the newly introduced concept of like-projection, and the M-stationary point is also investigated via a complementarity-type reformulation of the problem. The relationship between these  stationary points is discussed  as well. It turns out that all stationary points discussed in this work satisfy the necessary optimality conditions for the sparse bilinear least squares problem.
\vspace{3mm}

\noindent {\bf Key words:}\hspace{2mm} Bilinear least squares problem, sparsity constraint, first-order optimality condition, tangent cone, normal cone, like-projection
  \vspace{3mm}

\noindent {\bf AMS subject classifications:}\hspace{2mm} 90C30, 90C26, 90C46, 94A12

\end{abstract}

\section{Introduction}

Sparse optimization involves the restriction on the cardinality of the support of variables in either constraints or the objective  of an optimization problem (see, e.g., \cite{AP24, BBCS18, BS2018, BKS2016,TBLS24,Zhao2018}). Sparse optimization was largely motivated by the sparse signal recovery which has been extensively studied over the past decades (see, e.g., \cite{CT05,CRT06,D06,Elad2010,EK2012}). The problem for  sparse signal recovery can be cast as
\begin{eqnarray}\label{CS-basis}
\min_{x\in \mathbb{R}^n}\{\|Ax-b\|^2:\|x\|_0\leq s\},
\end{eqnarray} which can be referred to as the sparse least squares problem.
In this model, $A\in \mathbb{R}^{l\times n}$ is a matrix with $l\ll n$, $b\in R^l$ is a given vector (denoting the measurements of the signal to recover), and $s$ ($\ll n)$ is a given positive integer. Throughout the paper, $\|\cdot\|$ denotes the 2-norm of a vector and $\|\cdot\|_0$ is called the `$\ell_0$-norm' denoting  the number of nonzero components of a vector.
The model  (\ref{CS-basis}) is a special case of the following  optimization problem with a sparsity constraint:
\begin{eqnarray}\label{nonlinear-model}
\min\left\{f(x):\ \Vert x\Vert_0\le s\right\},
\end{eqnarray}
where $f:\mathbb{R}^n\to\mathbb{R}$ is a continuously differentiable function. The global sufficient and/or necessary optimality conditions for the problem  (\ref{CS-basis}) can be characterized in terms of such properties as the null space property  (see, e.g., \cite{CDD2008,FR2013}), restricted isometry property \cite{CT05}, range space propery of transposed matrices \cite{Zhao2013,Zhao2014}, and the mutual coherence and related properties \cite{Elad2010,Zhao2018}.  For instance, under the above-mentioned null or range space property, it is known that  the problem (\ref{CS-basis})  can be equivalent to an  $\ell_1$-norm-based convex optimization problem (e.g., the LASSO problem or $\ell_1$-minimization). This implies that the optimal solution of  (\ref{CS-basis}) can actually be  characterized via the optimality condition of its equivalent convex counterpart under such matrix properties.

Due to the generality of $f$ in (\ref{nonlinear-model}), it is generally more difficult to  characterize the optimality conditions for the solution of this problem. However, some recent progress has been made. Some work focuses on  developing the necessary optimality conditions  for (\ref{nonlinear-model}), and the others on certain sufficient optimality conditions for the problem. For instance, Beck and Eldar \cite{BE2013} investigated the  first-order necessary optimality conditions for the problem, including $L$-stationarity, basic feasibility and coordinate-wise minimality. They also developed certain numerical methods for (\ref{nonlinear-model}) to locate a point satisfying their optimality conditions. Pan et al. \cite{PXZ2015} studied the so-called  N-stationarity and T-stationarity conditions for (\ref{nonlinear-model}) and pointed out  that such stationarities are the necessary optimality conditions for this problem. They also developed a second-order sufficient optimality condition for (\ref{nonlinear-model}). Similar studies of optimality conditions for this problem can be found in \cite{MNS19,LS22,LS15} and numerical methods in \cite{ZXQ21,Z22,YLZ18,B13,BRB13,DS22,AP24,AL22}.

However, the model (\ref{nonlinear-model}) is far from being broad enough to cover all application scenarios. Motivated by structured signal recovery and other applications,  the optimization problems with more complex constraints than  (\ref{nonlinear-model})  have also been studied in the literature. Such problems can be formulated as
\begin{eqnarray}\label{general-model}
\begin{array}{cl}
\min & f(x)\\
\mbox{\rm s.t.}& g(x)\le0,\ h(x)=0,\\
 & \|x\|_0\le s,
\end{array}
\end{eqnarray}
where $g$ and $h$ are two continuously differentiable functions. Compared to (\ref{nonlinear-model}), it becomes more difficult to develop an optimality condition for the problem (\ref{general-model}). One way to achieve this is to borrow the analysis tool from the traditional optimization problems.  Burdakov, Kanzow and Schwartz \cite{BKS2016} proposed a complementarity-type reformulated model and thus established the M(Mordukhovich)-stationarity condition for  (\ref{general-model}). They also proposed a regularization method which turns out to be globally convergent to an M-stationary point of the problem. Based on such a  reformulation, Bucher and Schwartz \cite{BS2018} introduced  a second-order optimality condition that can  guarantee the local uniqueness of a M-stationary point. More studies on the optimality conditions of  (\ref{general-model}) can be found in \cite{LZ13,PXF2017,BBCS18,KRS21}.

In this paper, we study a sparse optimization problem which is also more general than (\ref{CS-basis}) but does not fall into the framework of (\ref{general-model}) . The model considered here is referred to as the sparse  bilinear least squares problem for which we develop certain necessary optimality conditions in this paper. Specifically, let $\mathscr{A}=(a_{ijk})\in\mathbb{R}^{l\times m\times n}$ be a given $3$-rd order tensor, $b\in\mathbb{R}^l$ be a given vector and $s,t$ be two fixed positive integers with $s<m$ and $t<n$. For $x\in\mathbb{R}^m$ and $ y\in\mathbb{R}^n,$ denote $\Gamma_1:={\rm supp(x)} = \{i: x_i \not =0\}$ and define $\gamma_x$ as the first index of $\Gamma_1$.  We consider the following model:
\begin{eqnarray}\label{originalmodel}
\begin{array}{cl}
\min & f(x,y)=\frac{1}{2}\|\mathscr{A} xy-b\|^2\\
\mbox{\rm s.t.} & e_{\gamma_x}^\top x=1,\\
& \|x\|_0\le s,\ \|y\|_0\le t,
\end{array}
\end{eqnarray}
where $e_{\gamma_x}\in\mathbb{R}^m$ is an $m$-dimensional identity column vector whose $\gamma_x$-th entry is $1$ and the others are $0.$ In this model,  $\mathscr{A}xy$ is short for $\mathscr{A}\times_2x\times_3y$, where  $\times_k$  ($k=2,3$) is the mod-$k$ product of $3$-rd order tensor. Precisely, for any $x\in\mathbb{R}^{m}$ and $y\in\mathbb{R}^{n}$, the entries of $\mathscr{A}\times_2x\in \mathbb{R}^{l\times n}$, $\mathscr{A}\times_3y\in \mathbb{R}^{l\times m}$ and $\mathscr{A}\times_2x\times_3y\in \mathbb{R}^{l}$ are given as
\begin{eqnarray*}
&& (\mathscr{A}\times_2x)_{ij}=\sum_{i_2\in \{1,2,\ldots,m\}}a_{ii_2j}x_{i_2},\ \forall i\in \{1,2,\ldots,l\}, \forall j\in \{1,2,\ldots,n\},\\
&& (\mathscr{A}\times_3y)_{ij}=\sum_{i_3\in \{1,2,\ldots,n\}}a_{iji_3}y_{i_3}, \ \forall i\in \{1,2,\ldots,l\}, \forall j\in \{1,2,\ldots,m\},\\
&& (\mathscr{A}\times_2x\times_3y)_{i}=\sum_{i_2\in \{1,2,\ldots,m\},i_3\in \{1,2,\ldots,n\}}a_{ii_2i_3}x_{i_2}y_{i_3}, \ \forall i\in \{1,2,\ldots,l\}.
\end{eqnarray*}
The constraint $e_{\gamma_x}^\top x=1$ is used to eliminate the multiplicity of the solution  to the problem and thus to guarantee the well-definedness of the  bilinear model \cite{EA2019,WZL2009}. In fact,  without the constraint $e_{\gamma_x}^\top x=1$, we see that for any $\alpha\in\mathbb{R}$, if $(\alpha x,y)$ solves the problem (\ref{originalmodel}), so does the point $(x,\alpha y)$.  Clearly, the well-definedness of the bilinear model can also ensured if $e_{\gamma_x}^\top x=1$ is replaced by $e_{\gamma_y}^\top y=1$, where $\gamma_y$ is the first index of ${\rm supp}(y). $ The constraint $e_{\gamma_x}^\top x=1$   makes  it more difficult to develop an optimality condition for the problem (\ref{originalmodel})  than  (\ref{nonlinear-model}) and (\ref{general-model}).

Problem (\ref{originalmodel})  arises in many scenarios, including blind deconvolution \cite{ARR2014,BPGD2014,BR2015,CMJ2017,LTR2018,CJ2019}, Hammerstein identification \cite{WZL2009,EA2019}, self-calibration \cite{GHLT1995,LLB2017} and matrix sensing \cite{AEB2006,LWB2018}. Let us take a look at two  examples.

\begin{example}
(Blind deconvolution) Generally, the mathematical model of blind deconvolution can be described as
\begin{eqnarray}\label{B-dec}
(Hx)\odot (Gy)=b,
\end{eqnarray}
where $H=(h_{ij})\in\mathbb{R}^{l\times m} $ and $ G=(g_{ik})\in\mathbb{R}^{l\times n}$ are two known dictionaries, $b\in\mathbb{R}^l$ is an observation vector, $x\in\mathbb{R}^m$ and $ y\in\mathbb{R}^n$ are two unknown vectors, and `$\odot$' denotes the Hadamard product (i.e., for any $u,v\in\mathbb{R}^l$, $w=u\odot v$ means $w_i=u_iv_i$ for  $i=1,\ldots,l$). In such scenarios as random mask imaging \cite{BR2015}, compressed sensing with unknown sensor gains \cite{CMJ2017, CJ2019} and  blind channel estimation \cite{LTR2018,ARR2014,LLB2017}, the sparse structure of $x,y$ needs to be considered separately. Thus the cardinality constraints $\Vert x\Vert_0\le s$ and $\Vert y\Vert_0\le t$ for certain integers $s<m$ and $t<n$ would come into play.  In this case,  (\ref{B-dec}) can be modeled as (\ref{originalmodel}) where $\mathscr{A}=(a_{ijk})\in\mathbb{R}^{l\times m\times n}$ with $a_{ijk}=g_{ij}h_{ik}$ for $i\in\{1,\ldots,l\}$, $j\in\{1,\ldots,m\}$ and $k\in\{1,\ldots,n\}$.
\end{example}
\begin{example}
(Matrix sensing)
Lee, Wu and Bresler \cite{LWB2018} solved rank-one matrix equations by using the sparse power factorization method when the row sparsity and column sparsity are both involved. Specifically, let $\mathcal{A}:\mathbb{R}^{m\times n}\to\mathbb{R}^{l}$ be a linear operator defined by $\mathcal{A}(X)=(\langle M_1,X\rangle,\ldots,\langle M_l,X\rangle)^\top$ for $X\in\mathbb{R}^{m\times n}$, where $M_i\in \mathbb{R}^{m\times n}$ are known measurements for $i\in\{1,\ldots,m\}$, and let  $b\in\mathbb{R}^l$ be the measurement vector and  $$S:=\left\{X\in\mathbb{R}^{m\times n}:\ X\ \text{has}\ s\ \text{nonzero rows and}\ t\ \text{nonzero columns}\right\}.$$ The sparse least squares problem for a matrix equation can be described as
\begin{eqnarray}\label{MEq}
\min_{X\in S}\ \Vert\mathcal{A}X-b\Vert^2.
\end{eqnarray}
Since $X$ is rank-one, there exist $\lambda\in\mathbb{R},x\in\mathbb{R}^m$ and $y\in\mathbb{R}^n$ such that $X=\lambda xy^\top$. Thus (\ref{MEq}) can be reformulated as the  sparse bilinear problem
\begin{eqnarray}\label{sBIP}
\min_{x\in S_1,y\in S_2}\ \Vert\lambda\mathcal{A}(xy^\top)-b\Vert^2,
\end{eqnarray}
where $S_1:=\left\{x\in\mathbb{R}^m:\ \Vert x\Vert_0\le s\right\}$ and $S_2:=\left\{y\in\mathbb{R}^n:\ \Vert y\Vert_0\le t\right\}$. So the problem above can be formulated  as (\ref{originalmodel}), where the tensor $\mathscr{A}=(a_{ijk})\in \mathbb{R}^{l\times m\times n}$ is defined by
$
a_{ijk}=(M_i)_{jk} $ for $ i\in \{1,\ldots,l\},  j\in \{1,\ldots,m\}$ and $ k\in \{1,\ldots,n\}.
$
\end{example}

The optimality conditions for (\ref{originalmodel}) have not been systematically studied at present.  This paper is devoted to this study  and the main contributions of this work include:
\begin{itemize}
\item The exact expressions of Bouligand (Clarke) tangent and normal cones for the problem (\ref{originalmodel}) at some feasible point are given. Based on this, we characterize the so-called $N^B\, (N^C)$-stationary point and $T^B\, (T^C)$-stationary point of  (\ref{originalmodel}). We show that an optimal solution of  (\ref{originalmodel}) must be an $N^B\, (N^C)$-stationary point and be an $T^B\, (T^C)$-stationary point for this problem.

\item The concept of coordinate-wise minima is given in accordance with the structure of the feasible set of (\ref{originalmodel}). We show that an optimal solution of (\ref{originalmodel}) is necessarily a CW-minimum.

\item Using the newly defined concept of  like-projection, L-like stationarity is introduced and characterized. In particular, we establish the first-order necessary condition for an optimal solution of (\ref{originalmodel}) to be  an L-like stationary point.

\item Based on a complementarity-type reformulation of (\ref{originalmodel}), we introduce and characterize the M-stationary point of  (\ref{originalmodel}). In particular, we also show that an optimal solution of  (\ref{originalmodel}) must be an M-stationary point.

\item We also discuss the relationship of the above-mentioned stationary points.
\end{itemize}

The rest of the paper is organized as follows. In Section \ref{sect2}, we introduce the notations, concepts and basic results that will be used in the subsequent analysis. In Section \ref{sect3}, we study the  $N^B (N^C)$-stationary points, $T^B (T^C)$-stationary points, coordinate-wise minima, L-like stationary points and M-stationary points  as well as  the related optimality conditions. The conclusions are given in Section \ref{sect4}.

\section{Preliminaries}\label{sect2}
\setcounter{equation}{0} \setcounter{assumption}{0}
\setcounter{theorem}{0} \setcounter{proposition}{0}
\setcounter{corollary}{0} \setcounter{lemma}{0}
\setcounter{definition}{0} \setcounter{remark}{0}
\setcounter{algorithm}{0}

Let $\mathbb{N}$ denote the set $\left\{1,2,\ldots,\right\}$. For fixed positive integers $m$ and $n$, we denote
\begin{eqnarray*}
I:=\left\{1,\ldots,m+n\right\}, I_1:=\left\{1,\ldots,m\right\}, I_2:=\left\{m+1,\ldots,m+n\right\}.
\end{eqnarray*}
Denote by  $e_{I}:=(1,\ldots,1)^\top\in\mathbb{R}^{m+n}$, $e_{I_1}:=(1,\ldots,1)^\top\in\mathbb{R}^m$ and $e_{I_2}:=(1,\ldots,1)^\top\in\mathbb{R}^n$. For any $i\in I$, let $e_i$ denote the $(m+n)$-dimensional column vector with the $i$-th entry being $1$ and the others being $0$. For fixed integer $l$, we denote $\mathscr{A}=(a_{i_1i_2i_3})\in\mathbb{R}^{l\times m\times n}$ as a $3$-rd order tensor with $i_1\in\left\{1,\ldots,l\right\}$, $i_2\in\left\{1,\ldots,m\right\}$ and $i_3\in\left\{1,\ldots,n\right\}$.

For any $x\in\mathbb{R}^m$,  $\Vert x\Vert$ denotes the $l_2$-norm of $x$,
$\Vert x\Vert_\infty$ denotes the infinite norm of $x$, ${\rm supp}(x):=\left\{i\in I_1:\ x_i\neq0\right\}$ denotes the support of $x$, and $ |x| $ denotes the absolute vector of $ x$, i.e.,  $|x|:=(|x_1|,\ldots,|x_m|)^\top$. We  use $M_i(x)$ to denote the value of the $i$-th largest entries of $x$ for any $i\in I_1$. For any $z:=(x^\top,y^\top)^\top\in\mathbb{R}^{m+n}$, we use $(x,y)$ to denote $(x^\top,y^\top)^\top$ for convenience,  and we use $\Gamma, \Gamma_1, \Gamma_2$ to denote respectively the following support sets:
\begin{eqnarray}\label{supphuang}
\Gamma:=\mbox{\rm supp}(z), ~ \Gamma_1:=\mbox{\rm supp}(x), ~ \Gamma_2:=\mbox{\rm supp}(y).
\end{eqnarray}
Let $\gamma_x$ be the index of the first nonzero entry of $x\in \mathbb{R}^m$, i.e., $\gamma_x:=\arg\min \left\{i\in I_1: x_i\neq 0 \right\}$.
Moreover, denote by
\begin{eqnarray}
F_1&:=&\left\{x\in\mathbb{R}^m:\ \Vert x\Vert_0\le s,\ e^\top_{\gamma_x}x=1\right\}, \label{fsbl}\\
F_2&:=&\left\{y\in\mathbb{R}^n:\ \Vert y\Vert_0\le t\right\}, \label{fsb2}\\
F&:=&F_1\times F_2=\left\{z=(x,y):\ x\in F_1,\ y\in F_2\right\}.\label{sprbl-set}
\end{eqnarray}
Clearly,  $F$ is the feasible set of (\ref{originalmodel}).

For a symmetric matrix $A\in\mathbb{R}^{m\times m}$, $\lambda_{max}(A)$ is the spectrum value of $A$, i.e.,
$$
\lambda_{max}(A):=\left\{\max|\lambda|:\ \lambda\in\mathbb{R}{\rm\ is\ an\ eigenvalue\ of\ }A\right\}.
$$
For an index set $\Lambda\subseteq I$, $|\Lambda|$ means the cardinality of $\Lambda$.
For a subset $S\subseteq\mathbb{R}^m$, let $\delta_S$ denote the indicator function of $S$. That is,
\begin{equation*}
\delta_S(x)=\left\{
\begin{array}{ll}
0 &\mbox{\rm if}\; x\in S,\\
\infty & \mbox{\rm if}\; x\notin S.
\end{array}
\right.
\end{equation*}

Tangent  and normal cones of the feasible set of optimization problems are important tools for the study of optimization problems (see,  e.g., \cite{C1990,RW1998}). The following concepts of tangent   and normal cones are defined in the sense of Bouligand and Clarke, respectively.
\begin{definition}\label{Bcone}
For any nonempty set $S\subseteq\mathbb{R}^n$, $d\in\mathbb{R}^n$ is called a Bouligand tangent vector of $S$ at $x\in S$ if there exist $\left\{x^k\right\}\subset S$ with ${\rm lim}_{k\to\infty}x^k=x$ and $\lambda_k\ge0$ for $k\in\mathbb{N}$ such that ${\rm lim}_{k\to\infty}\lambda_k(x^k-x)=d$. The set containing all Bouligand tangent vectors of $S$ at $x\in S$ is a cone, which is called the Bouligand tangent cone of $S$ at $x$, denote by $T^B_S(x)$. Moreover, the polar cone of $T^B_S(x)$ is called the Bouligand normal cone of $S$ at $x$, denote by $N^B_S(x)$, i.e.,
$$
N^B_S(x)=\left\{d\in\mathbb{R}^n:\ \left<d,u\right>\le0,\ \forall u\in T^B_S(x)\right\}.
$$
\end{definition}

\begin{definition}\label{Ccone}
For any nonempty set $S\subseteq\mathbb{R}^n$, $d\in\mathbb{R}^n$ is called a Clarke tangent vector of $S$ at $x\in S$ if for any $\left\{x^k\right\}\subset S$ with ${\rm lim}_{k\to\infty}x^k=x$ and $\left\{\lambda_k\right\}\subset\mathbb{R}_+$ with ${\rm lim}_{k\to\infty}\lambda_k=0$, there exists $\left\{y^k\right\}$ with ${\rm lim}_{k\to\infty}y^k=d$ such that $x^k+\lambda_ky^k\in S$ for any $k\in\mathbb{N}$. The set containing all Clarke tangent vectors of $S$ at $x\in S$ is called the Clarke tangent cone of $S$ at $x$, denote by $T^C_S(x)$. Correspondingly, the polar cone of the Clarke tangent cone of $S$ at $x$ is called the Clarke normal cone of $S$ at $x$, denote by $N^C_S(x)$, i.e,
$$
N^C_S(x)=\left\{d\in\mathbb{R}^n:\ \left<d,u\right>\le0,\ \forall u\in T^C_S(x)\right\}.
$$
\end{definition}

Throughout the remainder of the paper, we use $f: \mathbb{R}^{m+n}\rightarrow \mathbb{R}$ to denote the objective function of  (\ref{originalmodel}), $\nabla f(z)$ the gradient of $f$ at $z:=(x,y)\in\mathbb{R}^{m+n}$, and $\nabla_i f(z)$ the $i$-th component of $\nabla f(z)$ for $i\in I$. And for $S\subset I$, we use $\nabla_S f(z)$ to denote the sub-vector of $\nabla f(z) $ with only components $\nabla_i f(z)$ for $i\in S$.

\section{First-order optimality conditions}\label{sect3}
\setcounter{equation}{0} \setcounter{assumption}{0}
\setcounter{theorem}{0} \setcounter{proposition}{0}
\setcounter{corollary}{0} \setcounter{lemma}{0}
\setcounter{definition}{0} \setcounter{remark}{0}
\setcounter{algorithm}{0}

In this section, we investigate several classes of stationary points and the related optimality conditions for the problem (\ref{originalmodel}).

\subsection{N-type and T-type stationary points}
We first introduce the concepts of {\it N-type stationary point} and {\it T-type stationary point} via the tangent and normal cones described in Definitions \ref{Bcone} and \ref{Ccone}. Then we characterize these stationary points and establish the optimality conditions for problem  (\ref{originalmodel}).

\subsubsection{Basic concepts and lemma}

The definition of N-type stationary point of (\ref{originalmodel}) is an extension of the definition of critical point in \cite{ABS2013}.
\begin{definition}\label{def-N-stp}
Let $F$ be defined by (\ref{sprbl-set}). $z^*:=(x^*,y^*)\in\mathbb{R}^{m+n}$ is called an $N^\#$-stationary point of (\ref{originalmodel}) if
\begin{eqnarray*}
-\nabla f(z^*)\in N^\#_{F}(z^*),
\end{eqnarray*}
where $\#\in\left\{B,C\right\}$ means the Bouligand sense and Clarke sense, respectively.
\end{definition}

Next, the T-type stationary point of  (\ref{originalmodel}) is extended from \cite{CM1987}. 
\begin{definition}\label{def-T-stp}
Let $F$ be defined by (\ref{sprbl-set}). $z^*:=(x^*,y^*)\in\mathbb{R}^{m+n}$ is called a $T^\#$-stationary point of (\ref{originalmodel}) if
\begin{eqnarray*}
\Vert\nabla^\#_Ff(z^*)\Vert=0,
\end{eqnarray*}
where $\nabla^\#_Ff(z^*):={\rm argmin}\left\{\Vert\nabla f(z^*)+d\Vert:\ d\in T^\#_F(z^*)\right\}$ and $\#\in\left\{B,C\right\}$ means the Bouligand sense and Clarke sense, respectively.
\end{definition}

The following lemma will be frequently used in subsequent analyses.
\begin{lemma}\label{x-k-1}
For any $x\in F_1$ defined in (\ref{fsbl}) and $\left\{x^k\right\}\subset F_1$ satisfying ${\rm lim}_{k\to\infty}x^k=x$, there must exist $k_0\in \mathbb{N}$ such that the index for the first nonzero element of $x^k$ is equal to $\gamma_x$ for all $k>k_0$.
\end{lemma}

\noindent{\bf Proof}. Since the sequence $\left\{x^k\right\}$ converges to $x$, it follows that for any $\epsilon<1$, there exists $k_0\in \mathbb{N}$ such that $\Vert x^k-x\Vert_\infty<\epsilon$. Choose $k>k_0$ arbitrarily and denote the first element of ${\rm supp}(x^k)$ by $i^*$.
\begin{itemize}
\item If $i^*<\gamma_x$, then we have $x^k_{i^*}=1$ and $x_{i^*}=0$ by the definition of $F_1$ since $x\in F_1$ and $\left\{x^k\right\}\subset F_1$. Thus $\Vert x^k-x\Vert_\infty\ge|x^k_{i^*}-x_{i^*}|=1>\epsilon$, which contradicts the result $\Vert x^k-x\Vert_\infty<\epsilon$.
\item If $i^*>\gamma_x$, then we have $x^k_{\gamma_x}=0$ and $x_{\gamma_x}=1$, and hence, $\Vert x^k-x\Vert_\infty\ge|x^k_{\gamma_x}-x_{\gamma_x}|=1>\epsilon$, which also leads to a contradiction.
\end{itemize}
Thus, $i^*=\gamma_x$. By the arbitrariness of $k$, we conclude that the index for the first nonzero element of $x^k$ is equal to $\gamma_x$ for any $k>k_0$.
\ep

\subsubsection{ First-order optimality conditions by  Bouligand tangent and normal cones}
Firstly, we describe the precise forms of the tangent and normal cones of $F$ at a point in $F$ in the sense of Bouligand. The following two lemmas are useful for proving the desired result.

\begin{lemma}(\cite[Proposition 6.41]{RW1998})\label{B-product}
Assume that $C_1\subseteq\mathbb{R}^m$ and $C_2\subseteq\mathbb{R}^n$, then for any $z=(x,y)\in C_1\times C_2\subseteq\mathbb{R}^{m+n}$, we have
\begin{eqnarray*}
T_{C_1\times C_2}^B(z)=T_{C_1}^B(x)\times T_{C_2}^B(y)\quad \mbox{\rm and}\quad N_{C_1\times C_2}^B(z)=N_{C_1}^B(x)\times N_{C_2}^B(y).
\end{eqnarray*}
\end{lemma}

The above lemma indicates that the Bouligand tangent and normal cones are separable under the Cartesian product of sets. In fact, denote $C:=C_1\times C_2$, by the definition of Bouligand tangent cone at some point $z=(x,y)\in C$, we have
\begin{eqnarray*}
T^B_C(z)=\left\{d\in\mathbb{R}^{m+n}:\ d=\lim_{k\to\infty}\lambda_k(z^k-z)\ \mbox{\rm where}\ z^k\overset{C}\to z,\left\{\lambda_k\right\}\ge0\right\}.
\end{eqnarray*}
Since $C=C_1\times C_2$,  we have $d=(d_1,d_2)\in T^B_C(z)$ if and only if for $\left\{\lambda_k\right\}\ge0$,
\begin{eqnarray*}
d_1=\lim_{k\to\infty}\lambda_k(x^k-x)\ \mbox{\rm where}\ x^k\overset{C_1}\to x\quad \mbox{\rm and}\quad d_2=\lim_{k\to\infty}\lambda_k(y^k-y)\ \mbox{\rm where}\ y^k\overset{C_2}\to y.
\end{eqnarray*}
That is $d_1\in T^B_{C_1}(x)$ and $d_2\in T^B_{C_2}(y)$. The result for the normal cone follows directly from the polarity relation. This lemma provides an approach to characterize the Bouligand tangent and normal cones of the set $F$, that is, by first characterizing those of $F_1$ and $F_2$ separately and then taking their Cartesian product. Fortunately, the characterization of $F_2$ has already been established in the existing literature, which is given in the following lemma.

\begin{lemma}(\cite[Theorem 2.1]{PXZ2015})\label{BNT-F2}
For any $y\in F_2$ with $F_2$ being defined by (\ref{fsb2}), the Bouligand tangent cone and corresponding normal cone of $F_2$ at $y$ are
\begin{eqnarray}
T^{B}_{F_2}(y)&=&\left\{d\in\mathbb{R}^n:\ \Vert d\Vert_0\le t,\ \Vert y+\mu d\Vert_0\le t,\ \forall\mu\in\mathbb{R}\right\}\label{btc-y-1}\\
&=&\bigcup_{\Lambda}{\rm span}\left\{e_i,i\in\Lambda\supseteq\Gamma_2,|\Lambda|\le t\right\},\label{btc-y-2}
\end{eqnarray}
\begin{equation}\label{bnc-y}
N^B_{F_2}(y)=\left\{
\begin{array}{lcl}
\left\{d\in\mathbb{R}^{n}:\ d_i=0,\ i\in\Gamma_2\right\}&&\mbox{\rm if}\ |\Gamma_2|=t,\\
\left\{0\right\}&&\mbox{\rm if}\ |\Gamma_2|<t,
\end{array}
\right.
\end{equation}
respectively.
\end{lemma}

\begin{theorem}\label{BNT-F1}
For any $x\in F_1$ with $F_1$ being defined by (\ref{fsbl}), then the Bouligand tangent cone and normal cone of $F_1$ at $x$ can be written as
\begin{eqnarray}
T^B_{F_1}(x)&=&\left\{d\in\mathbb{R}^m:\ \Vert d\Vert_0\le s-1;\ d_j=0,\ \forall j\le\gamma_x;\ x+\mu d\in F_1,\  \forall\mu\in\mathbb{R}\right\}\label{btc-f1}\\
&=&\bigcup_\Lambda{\rm span}\left\{e_i:\ i\in\Lambda\supseteq\Gamma_1\setminus\left\{\gamma_x\right\};\ j\notin\Lambda,\ \forall j\le\gamma_x;\ |\Lambda|\le s-1\right\},\label{btc-f1-2}
\end{eqnarray}
\begin{equation}\label{bnc-f1}
N^{B}_{F_1}(x)=\left\{
\begin{array}{lcl}
\left\{d\in\mathbb{R}^{m}:\ d_i=0,\ i\in\Gamma_1\setminus\left\{\gamma_x\right\}\right\}&&\mbox{\rm if}\ |\Gamma_1|=s,\\
\left\{d\in\mathbb{R}^{m}:\ d_i=0,\ i>\gamma_x\right\}&&\mbox{\rm if}\ |\Gamma_1|<s,
\end{array}
\right.
\end{equation}
respectively.
\end{theorem}
\noindent{\bf Proof}. Let $ D$ denote the set on the right-hand side of (\ref{btc-f1}). Thus it is sufficient to show that $T^B_{F_1}(x)=D$. On one hand, for any $d\in T^B_{F_1}(x)$, it follows from Definition \ref{Bcone} that there exist $\left\{x^k\right\}\subset F_1$ satisfying $\lim_{k\to\infty}x^k=x$ and $\lambda_k\geq 0$ for $k\in \mathbb{N}$ such that $\lim_{k\to\infty}\lambda_k(x^k-x)=d$. Then there exists $k_0\in \mathbb{N}$ such that ${\rm supp}(x)\subseteq{\rm supp}(x^k)$ and ${\rm supp}(d)\subseteq{\rm supp}(x^k-x)\subseteq{\rm supp}(x^k)$ for any $k>k_0$, and hence,
 \begin{equation}
{\rm supp}(d)\subseteq{\rm supp}(x^k),\ {\rm supp}(d)\cup{\rm supp}(x)\subseteq{\rm supp}(x^k). \label{suppd1}
\end{equation}
These imply that
\begin{eqnarray}\label{zk-mud-F1}
\Vert x+\mu d\Vert_0\le s,\ \forall\mu\in\mathbb{R}.
\end{eqnarray}
In addition, by Lemma \ref{x-k-1}, there exists $k_1\in\mathbb{N}$ such that $x^k_{\gamma_x}=1$ for any $k>k_1$. Thus,
\begin{eqnarray}\label{d-gx-0}
d_{\gamma_x}={\rm lim}_{k\to\infty}\lambda_k(x^k_{\gamma_x}-x_{\gamma_x})=0.
\end{eqnarray}
Note that  for any $j\in I_1$ with $j<\gamma_x$,
\begin{eqnarray}\label{d-gxm-0}
d_j={\rm lim}_{k\to\infty}\lambda_k(x^k_j-x_j)={\rm lim}_{k\to\infty}\lambda_k(0-0)=0.
\end{eqnarray}
Thus, for any $\mu\in\mathbb{R},$ one has
\begin{eqnarray}
(x+\mu d)_{\gamma_x}&=&x_{\gamma_x}+\mu d_{\gamma_x}=1,\label{z-mud-gx-1}\\
(x+\mu d)_j&=&x_j+\mu d_j=0,\ \forall j<\gamma_x.\label{z-mud-gx-0}
\end{eqnarray}
Combining (\ref{suppd1})-(\ref{z-mud-gx-0}), we can obtain that $z+\mu d\in F$ for any $\mu\in\mathbb{R}$, $\Vert d\Vert_0\le s-1$ and $d_j=0$ for any $j\in I_1$ with $j<\gamma_x$. This means $d\in D$, and thus $T^B_{F_1}(x)\subseteq D$. On the other hand, for any $d\in D$, we take $\left\{\lambda_k\right\}$ such that $\lim_{k\to\infty}\lambda_k=\infty$, and define
$$
x^k:=x+\frac{1}{\lambda_k}d,\ \forall k\in\mathbb{N}.
$$
Then, we have $x^k\in F_1$ for any $k\in\mathbb{N}$, $\lim_{k\to\infty}x^k=x$ and $d=\lim_{k\to\infty}\lambda_k(x^k-x)$. It follows  that $d\in T^B_{F_1}(x)$. Thus $D\supseteq T^B_{F_1}(x)$. So (\ref{btc-f1}) holds, we now further show that (\ref{btc-f1-2}) holds.

Let $P$ denote the set on the right-hand side of (\ref{btc-f1-2}). We now prove (\ref{btc-f1-2}), i.e., $D=P$. On one hand, for any $d\in D$, we define
$
\Lambda^d:={\rm supp}(d)\cup\left(\Gamma_1\setminus\left\{\gamma_x\right\}\right).
$
It is easy to see that
$$
d\in{\rm span}\left\{e_i:\ i\in\Lambda^d\right\}, ~ \Gamma_1\setminus\left\{\gamma_x\right\}\subseteq\Lambda^d, ~ \gamma_x\notin\Lambda^d
$$
and $j\notin\Lambda^d$ for any $j\in I_1$ with $j<\gamma_x$. In addition, since $\Vert x+\mu d\Vert_0\le s$ for any $\mu\in\mathbb{R}$, one has
\begin{eqnarray*}
|{\rm supp}(d)\cup{\rm supp}(x)|\le s,
\end{eqnarray*}
which together with $\gamma_x\in{\rm supp}(x)$ imply that $|\Lambda^d|\le s-1$. Thus it follows from (\ref{btc-f1-2}) that ${\rm span}\left\{e_i:\ i\in\Lambda^d\right\}\subseteq P$, and thus   $d\in P$. So  $D\subseteq P$. On the other hand, for any $d\in P$, it follows from (\ref{btc-f1-2}) that there exists $\bar{\Lambda}^d$ with $|\bar{\Lambda}^d|\le s-1$, $\Gamma_1\setminus\left\{\gamma_x\right\}\subseteq\bar{\Lambda}^d$ and $j\notin\bar{\Lambda}^d$ for any $j\in I_1$ with $j<\gamma_x$ such that ${\rm supp}(d)\subseteq\bar{\Lambda}^d$. Therefore, $|{\rm supp}(d)|\le s-1$ and $d_j=0$ for any $j\in I_1$ with $j<\gamma_x$. Moreover, ${\rm supp}(x+\mu d)\subseteq\Gamma_1\cup{\rm supp}(d)\subseteq\bar{\Lambda}^d\cup\left\{\gamma_x\right\}$ for any $\mu\in\mathbb{R}$. Thus, $\Vert x+\mu d\Vert_0\le s$, $(x+\mu d)_{\gamma_x}=1$ and $(x+\mu d)_j=0$ for any $j\in I_1$ with $j<\gamma_x$. That is, $x+\mu d\in F_1$ for any $\mu\in\mathbb{R}$. This, together with  $|{\rm supp}(d)|\le s-1$ and $ d_j=0$ for any $j\in I_1$ with $j<\gamma_x$, implies that $d\in D$. Thus, $P\subseteq D$. So $D=P$, and thus  (\ref{btc-f1-2}) holds.

Now, we show  (\ref{bnc-f1}), which holds trivially when $s=1$. Specifically, one must have $\Gamma_1=\left\{\gamma_x\right\}$ for any $x\in F_1$ with $s=1$. Then at this point $x$, the tangent cone $T^B_{F_1}(x)$ reduces to $\left\{0\right\}$. It follows from the definition of normal cone that $N^B_{F_1}(x)=\mathbb{R}^m$, which implies that (\ref{bnc-f1}) holds due to $\Gamma_1=\left\{\gamma_x\right\}$. Now, we turn to the case $s>1$. In this situation, (\ref{bnc-f1}) can be verified by considering two separate cases.

  \emph{Case 1}: $|\Gamma_1|=s$.
On one hand, for any $h\in T^B_{F_1}(x)$, it follows from (\ref{btc-f1}) that $h_j=0$ for any $j\notin\Gamma_1\setminus\left\{\gamma_x\right\}$. Thus, for any $d\in\left\{d\in\mathbb{R}^{m}:\ d_i=0,\forall i\in\Gamma_1\setminus\left\{\gamma_x\right\}\right\}$, we have $\left<d,h\right>=0\le0$, and hence $d\in N^B_{F_1}(x)$ by the definition of $N^B_{F_1}(x)$. So in this case
$$
\left\{d\in\mathbb{R}^{m}:\ d_i=0,\forall i\in\Gamma_1\setminus\left\{\gamma_x\right\}\right\}\subseteq N^B_{F_1}(x).
$$
On the other hand, for any $d\in N^B_{F_1}(x)$, it follows from (\ref{btc-f1}) that $d\in D$. If $$d\notin\left\{d\in\mathbb{R}^{m}:\ d_i=0, \forall i\in\Gamma_1\setminus\left\{\gamma_x\right\}\right\},$$  then there exists $i_0\in\Gamma_1\setminus\left\{\gamma_x\right\}$ such that $d_{i_0}\neq0$. Define $h\in\mathbb{R}^{m}$ such that $h_{i_0}=d_{i_0}$ and other entries are zero. Then, it is easy to see that
 $\Vert h\Vert_0\le s-1$ and $x+\mu h\in F_1$ due to $i_0\in\Gamma_1\setminus\left\{\gamma_x\right\}$, and that
 $h_j=0$ for all $j\in I_1$ with $j\le\gamma_x$ by the definition of $\gamma_x. $
 Thus, $h\in D$, i.e., $h\in T^B_{F_1}(x)$. However, $\left<d,h\right>=d_{i_0}^2>0$, which contradicts $d\in N^B_{F_1}(x)$ and $h\in T^B_{F_1}(x)$. Thus, $d\in\left\{d\in\mathbb{R}^{m}:\ d_i=0,\forall i\in\Gamma_1\setminus\left\{\gamma_x\right\}\right\}$. So
$$N^B_{F_1}(z)\subseteq \left\{d\in\mathbb{R}^{m}:\ d_i=0,\forall i\in\Gamma_1\setminus\left\{\gamma_x\right\}\right\}.
$$
Therefore,  (\ref{bnc-f1}) holds in this case.

\emph{Case 2}: $|\Gamma_1|<s$.
On one hand, for any
$$
d\in\hat{D}:=\left\{d\in\mathbb{R}^{m}:\ d_i=0,\forall i>\gamma_x\right\},
$$
one has $d_i=0$ for any $i\in \left\{i\in I_1:i>\gamma_x\right\}$. Since for any $h\in T^B_{F_1}(x)$, it follows from (\ref{btc-f1}) that $h_j=0$ for all $j\in I_1$ with $j\le\gamma_x$, we have $\left<d,h\right>=0\le0$ for any $h\in T^B_{F_1}(x)$. Then we have $d\in N^B_{F_1}(x)$. So $\hat{D}\subseteq N^B_{F_1}(x)$. On the other hand, suppose there exists $d\in N^B_{F_1}(x)$ but $d\notin\hat{D}$. Then there exists $i_0\in \left\{i\in I_1:i>\gamma_x\right\}$ such that $d_{i_0}\neq0$. We define $h$ same as Case 1. Similar to the analysis in Case 1, we obtain a contradiction as well. Thus, $d\in\hat{D}$. So $\hat{D}\supseteq N^B_{F_1}(x)$. Therefore, $\hat{D}=N^B_{F_1}(x)$. It means (\ref{bnc-f1}) still holds in this case. The proof of the theorem is complete. \ep

\begin{theorem}
For any $z=(x,y)\in F$ with $F$ being defined by (\ref{sprbl-set}), the Bouligand tangent cone and corresponding normal cone of $F$ at $z$ are given by
\begin{align}
T^B_F(z)&=\left\{
\begin{array}{ll}
d=\left(d_1,d_2\right)\in\mathbb{R}^{m+n}:&\Vert d_1\Vert_0\le s-1,\ \Vert d_2\Vert_0\le t;\\
 & z+\mu d\in F, \ d_j=0 \textrm{ for all } j\le\gamma_x
\end{array}
\right\} \label{btc-z-1}\\
&=\bigcup_{\Lambda}{\rm span}\left\{
\begin{array}{ll}
\{e_i:i\in\Lambda\}:& \Gamma\setminus\left\{\gamma_x\right\}\subset \Lambda, j\notin\Lambda\textrm{ for all } j\le\gamma_x,\\
& |\Lambda\cap I_1|\le s-1,\ |\Lambda\cap I_2|\le t\\
\end{array}
\right\}, \label{btc-z-2}
\end{align}

\begin{equation}  N^B_F(z)=\left\{
\begin{array}{ll}
\left\{d\in\mathbb{R}^{m+n}:\ d_i=0,\forall i\in\Gamma\setminus\left\{\gamma_x\right\}\right\}&\mbox{\rm if}\ |\Gamma_1|=s\ \mbox{\rm and}\ |\Gamma_2|=t,\\
\left\{d\in\mathbb{R}^{m+n}:\ d_i=0,\forall i\in\left(\Gamma_1\setminus\left\{\gamma_x\right\}\right)\cup I_2\right\} &\mbox{\rm if}\ |\Gamma_1|=s\ \mbox{\rm and}\ |\Gamma_2|<t,\\
\left\{d\in\mathbb{R}^{m+n}:\ d_i=0,\forall i\in\left\{\gamma_x+1,\ldots,m\right\}\cup \Gamma_2\right\}&\mbox{\rm if}\ |\Gamma_1|<s\ \mbox{\rm and}\ |\Gamma_2|=t,\\
\left\{d\in\mathbb{R}^{m+n}:\ d_i=0,\forall i\in\left\{i\in I:i>\gamma_x\right\}\right\} &\mbox{\rm if}\ |\Gamma_1|<s\ \mbox{\rm and}\ |\Gamma_2|<t,
\end{array}
\right.  \label{bnc-z}
\end{equation}
respectively, where $\Gamma,\Gamma_1$ and $\Gamma_2$ are defined by (\ref{supphuang}).
\end{theorem}

The results are directly from Lemma \ref{B-product}, Lemma \ref{BNT-F2} and Theorem \ref{BNT-F1}. Next, we present a sufficient and necessary condition  for ${\rm N^B}$- and ${\rm T^B}$-stationary points, respectively.  The first one below can be obtained directly from (\ref{bnc-z}).

\begin{theorem}\label{thm-N-B}
Let $z^*=(x^*,y^*)\in F,$ where $F$ is defined by (\ref{sprbl-set}). Then
$z^*$ is an ${\rm N^B}$-stationary point of (\ref{originalmodel}) if and only if $\nabla f(z^*)\in \mathbb{R}^{m+n}$ with
\begin{eqnarray}\label{N-B-c1}
\left\{\begin{array}{ll}
\nabla_if(z^*)=0,\forall i\in\Gamma^*\setminus\left\{\gamma_{x^*}\right\}&\ \mbox{\rm if}\; |\Gamma^*_1|=s\; \mbox{\rm and}\; |\Gamma^*_2|=t, \\
\nabla_if(z^*)=0,\forall i\in(\Gamma^*_1\setminus\left\{\gamma_{x^*}\right\})\cup I_2&\ \mbox{\rm if}\; |\Gamma^*_1|=s\; \mbox{\rm and}\; |\Gamma^*_2|<t, \\
\nabla_if(z^*)=0,\forall i\in\left\{\gamma_{x^*}+1,\ldots,m\right\}\cup\Gamma^*_2&\ \mbox{\rm if}\; |\Gamma^*_1|<s\; \mbox{\rm and}\; |\Gamma^*_2|=t,\\
\nabla_if(z^*)=0,\forall i\in\left\{i\in I:i>\gamma_{x^*}\right\}&\ \mbox{\rm if}\; |\Gamma^*_1|<s\; \mbox{\rm and}\; |\Gamma^*_2|<t,
\end{array}\right.
\end{eqnarray}
where $\Gamma^*,\Gamma^*_1,\Gamma^*_2$ are defined by
\begin{eqnarray}\label{gm-star}
\Gamma^*:=\mbox{\rm supp}(z^*), ~ \Gamma^*_1:=\mbox{\rm supp}(x^*), ~\Gamma^*_2:=\mbox{\rm supp}(y^*).
\end{eqnarray}
\end{theorem}

\begin{theorem}\label{thm-T-B}
Let $z^*=(x^*,y^*)\in F,$ where $F$ is defined by (\ref{sprbl-set}). Then $z^*$ is a ${\rm T^B}$-stationary point if and only if (\ref{N-B-c1}) holds.
\end{theorem}

\noindent{\bf Proof}. Define $\Gamma^*,\Gamma^*_1,\Gamma^*_2$ as (\ref{gm-star}). We consider the following four cases.

\emph{Case 1}: $|\Gamma^*_1|=s$ and $|\Gamma^*_2|=t$. For this case, we prove that $\nabla^B_Ff(z^*)=0$ if and only if $\nabla_if(z^*)=0$ for all $i\in\Gamma^*\setminus\left\{\gamma_{x^*}\right\}$.

($\Rightarrow$) Suppose there exists $i_0\in\Gamma^*\setminus\left\{\gamma_{x^*}\right\}$ such that $\nabla_{i_0}f(z^*)\neq0$. Pick $\bar{d}=(\bar{d}_1,\bar{d}_2)\in\mathbb{R}^{m+n}$ such that $\bar{d}_{i_0}=-\nabla_{i_0}f(z^*)$ and other entries of $\bar{d}$ are zero. Since $\Vert \bar{d}_1\Vert_0\le s-1$ (note that $d$ is well-defined here: when $s=1$, since $i_0\in\Gamma^*\setminus\left\{\gamma_x\right\}$, we have $d_1=0$ i.e., $\Vert d_1\Vert_0= 0$), $\Vert\bar{d}_2\Vert_0\le t$, $z^*+\mu\bar{d}\in F$ for all $\mu\in\mathbb{R}$ and $\bar{d}_j=0$ for all $j\in I_1$ with $j\le\gamma_{x^*}$, it follows from (\ref{btc-z-1}) that $\bar{d}\in T^B_F(z^*)$. Then, by Definition \ref{def-T-stp}, one  has
\begin{eqnarray}\label{argmin-nablaF}
\Vert\nabla^B_Ff(z^*)+\nabla f(z^*)\Vert^2\le\Vert\bar{d}+\nabla f(z^*)\Vert^2.
\end{eqnarray}
Since $\nabla^B_Ff(z^*)=0$, (\ref{argmin-nablaF}) is equivalent to
\begin{eqnarray}\label{i0le}
\sum_{1\le i \le m+n}(\nabla_if(z^*))^2\le \sum_{1\le i \le m+n,i\neq i_0}(\nabla_if(z^*))^2.
\end{eqnarray}
This implies that $(\nabla_{i_0}f(z^*))^2\le0$,  a contradiction.

($\Leftarrow$) Suppose there exists $\bar{d}\in T^B_F(z^*)$ and $\bar{d}\neq0$ such that
\begin{eqnarray}\label{bar-d-le}
\Vert\bar{d}+\nabla f(z^*)\Vert^2\le\Vert0+\nabla f(z^*)\Vert^2.
\end{eqnarray}
It follows from (\ref{btc-z-1}) that ${\rm supp}(\bar{d})\subseteq\Gamma^*\setminus\left\{\gamma_{x^*}\right\}$. Thus (\ref{bar-d-le}) is equivalent to
$$
\sum_{i\in\Gamma^*\setminus\left\{\gamma_{x^*}\right\}}(\bar{d}_i+\nabla_if(z^*))^2\le\sum_{i\in\Gamma^*\setminus\left\{\gamma_{x^*}\right\}}(\nabla_if(z^*))^2.
$$
Since $\nabla_if(z^*)=0$ for all $i\in\Gamma^*\setminus\left\{\gamma_{x^*}\right\}$, the  inequality above implies that $\sum_{i\in\Gamma^*\setminus\left\{\gamma_{x^*}\right\}}(\bar{d}_i)^2\le0$, which leads to a contradiction.

\emph{Case 2}: $|\Gamma^*_1|=s$ and $|\Gamma^*_2|<t$. In this case, we prove that $\nabla^B_Ff(z^*)=0$ if and only if $\nabla_if(z^*)=0$ for all $i\in(\Gamma^*_1\setminus\left\{\gamma_{x^*}\right\})\cup I_2$.

($\Rightarrow$) Suppose there exists $i_0\in(\Gamma^*_1\setminus\left\{\gamma_{x^*}\right\})\cup I_2$ such that $\nabla_{i_0}f(z^*)\neq0$. Pick $\bar{d}=(\bar{d}_1,\bar{d}_2)\in\mathbb{R}^{m+n}$ same as Case 1, we also have $\bar{d}\in T^B_F(z^*)$. Then, by Definition \ref{def-T-stp}, the inequality (\ref{argmin-nablaF}) holds.
Since $\nabla^B_Ff(z^*)=0$, (\ref{argmin-nablaF}) is equivalent to (\ref{i0le}), which implies that $(\nabla_{i_0}f(z^*))^2\le0$. This is a contradiction.

($\Leftarrow$) Suppose there exists $\bar{d}\in T^B_F(z^*)$ and $\bar{d}\neq0$ such that (\ref{bar-d-le}) holds.
It follows from (\ref{btc-z-1}) that ${\rm supp}(\bar{d})\subseteq(\Gamma^*_1\setminus\left\{\gamma_{x^*}\right\})\cup I_2$. Thus  (\ref{bar-d-le}) is equivalent to
$$
\sum_{i\in(\Gamma^*_1\setminus\left\{\gamma_{x^*}\right\})\cup I_2}(\bar{d}_i+\nabla_if(z^*))^2\le\sum_{i\in(\Gamma^*_1\setminus\left\{\gamma_{x^*}\right\})\cup I_2}(\nabla_if(z^*))^2,
$$
which implies $\sum_{i\in(\Gamma^*_1\setminus\left\{\gamma_{x^*}\right\})\cup I_2}(\bar{d}_i)^2\le0$, a contradiction.

\emph{Case 3}: $|\Gamma^*_1|<s$ and $|\Gamma^*_2|=t$. In this case, we prove that $\nabla^B_F(z^*)=0$ if and only if $\nabla_if(z^*)=0$ for all $i\in\left\{\gamma_{x^*}+1,\ldots,m\right\}\cup\Gamma^*_2$.

($\Rightarrow$) Suppose there exists $i_0\in\left\{\gamma_{x^*}+1,\ldots,m\right\}\cup\Gamma^*_2$ such that $\nabla_{i_0}f(z^*)\neq0$. Pick $\bar{d}=(\bar{d}_1,\bar{d}_2)\in\mathbb{R}^{m+n}$ same as Case 1, we also have $\bar{d}\in T^B_F(z^*)$. Then, by Definition \ref{def-T-stp}, one must have (\ref{argmin-nablaF}) holds. Similar to the analysis in Case 1, it follows from $\nabla^B_Ff(z^*)=0$ that $(\nabla_{i_0}f(z^*))^2\le0$. This is a contradiction.

($\Leftarrow$) Suppose there exists $\bar{d}\in T^B_F(z^*)$ and $\bar{d}\neq0$ such that (\ref{bar-d-le}) holds.
It follows from (\ref{btc-z-1}) that ${\rm supp}(\bar{d})\subseteq\left\{\gamma_{x^*}+1,\ldots,m+n\right\}$. Thus (\ref{bar-d-le}) implies $\sum_{i=\gamma_{x^*}+1}^{m+n}(\bar{d}_i)^2\le0$, a contradiction.

\emph{Case 4}: $|\Gamma^*_1|<s$ and $\Gamma^*_2<t$. In this case, we prove that $\nabla^B_F(z^*)=0$ if and only if $\nabla_if(z^*)=0$ for all $i\in \left\{i\in I:i>\gamma_{x^*}\right\}$.

($\Rightarrow$) Suppose there exists $i_0\in\left\{i\in I:i>\gamma_{x^*}\right\}$ such that $\nabla_{i_0}f(z^*)\neq0$. Pick $\bar{d}=(\bar{d}_1,\bar{d}_2)\in\mathbb{R}^{m+n}$ same as Case 1, we also have $\bar{d}\in T^B_F(z^*)$. Then, by the same analysis as in Case 1, one can also obtain that $(\nabla_{i_0}f(z^*))^2\le0$, a contradiction.

($\Leftarrow$) Suppose there exists $\bar{d}\in T^B_F(z^*)$ and $\bar{d}\neq0$ such that (\ref{bar-d-le}) holds.
It follows from (\ref{btc-z-1}) that ${\rm supp}(\bar{d})\subseteq\left\{i\in I: i>\gamma_{x^*}\right\}$. Thus (\ref{bar-d-le}) implies $\sum_{i\in\left\{i\in I:i>\gamma_{x^*}\right\}}(\bar{d}_i)^2\le0$,  a contradiction again. Therefore the proof of the theorem is complete. \ep

According to Theorems \ref{thm-N-B} and \ref{thm-T-B}, the next result follows immediately.

\begin{corollary}\label{thm-huang}
$z^*\in\mathbb{R}^{m+n}$ is an ${\rm N^B}$-stationary point of (\ref{originalmodel}) if and only if $z^*$ is a ${\rm T^B}$-stationary point of (\ref{originalmodel}).
\end{corollary}

Finally, we discuss the relation between ${\rm N^B}$-stationary point (${\rm T^B}$-stationary point) of (\ref{originalmodel}) and the optimal solution of (\ref{originalmodel}).

\begin{theorem}\label{B-nec}
If $z^*=(x^*,y^*)$ is an optimal solution of (\ref{originalmodel}), then $z^*$ is an ${\rm N^B}$-stationary point, and $z^*$ is also a ${\rm T^B}$-stationary point.
\end{theorem}

\noindent{\bf Proof}. By Corollary \ref{thm-huang}, we only need to show that $z^*$ is an ${\rm N^B}$-stationary point. We will use Theorem \ref{thm-N-B} to show this result. Since $z^*=(x^*,y^*)$ is an optimal solution of (\ref{originalmodel}), it follows that $z^*\in F$. Define $\Gamma^*,\Gamma^*_1,\Gamma^*_2$ by (\ref{gm-star}). We now show the desired result by using (\ref{N-B-c1}).

Let $g_i(u):=f(z^*+ue_i)$ for any $i\in I$ and $u\in \mathbb{R}$, where   $f$ is defined by (\ref{originalmodel}). Then, for any $i\in I$, we have
\begin{eqnarray*}
g_i(u)=\left\{\begin{array}{ll}
\Vert\mathscr{A}(x^*+ue_i)y^*-b\Vert^2 & \mbox{\rm if}\; i\in I_1,\\
\Vert\mathscr{A}x^*(y^*+ue_i)-b\Vert^2 & \mbox{\rm if}\; i\in I_2,
\end{array}\right.
\end{eqnarray*}
which  is a strictly convex function for any $i\in I$. Thus $g_i(\cdot)$ has a unique minimizer that satisfies $g'_i(u)=0$. In what follows, we divide the proof into four cases.

 \emph{Case 1}: $\Vert x^*\Vert_0=s$ and $\Vert y^*\Vert_0=t$. For any $i\in\Gamma^*\setminus\left\{\gamma_{x^*}\right\}$, we have
$
{\rm argmin}\left\{g_i(u):\ u\in\mathbb{R}\right\}=0.
$
Otherwise, there exists $u_0\neq0$ such that $z^*+u_0e_i\in F$ and $f(z^*+u_0e_i)<f(z^*)$ which yields a contradiction. Thus $g'_i(0)=0$ for any $i\in\Gamma^*\setminus\left\{\gamma_{x^*}\right\}$. Moreover,
$$
g'_i(0)={\rm lim}_{u\to0}\frac{f(z^*+ue_i)-f(z^*)}{u}=\nabla_if(z^*),\ \forall i\in I.
$$
Thus $\nabla_if(z^*)=g'_i(0)=0$ for any $i\in\Gamma^*\setminus\left\{\gamma_{x^*}\right\}$. This, together with Theorem \ref{thm-N-B}, implies that $z^*$ is an ${\rm N^B}$-stationary point when $\Vert x^*\Vert_0=s$ and $\Vert y^*\Vert_0=t$.

\emph{Case 2}: $\Vert x^*\Vert_0=s$ and $\Vert y^*\Vert_0<t$. For any $i\in\left(\Gamma^*_1\setminus\left\{\gamma_{x^*}\right\}\right)\cup I_2$, we have $z^*+ue_i\in F$ for any $u\in\mathbb{R}$. If $g'_i(0)\neq0$ for any $i\in\Gamma^*_1\cup I_2\setminus\left\{\gamma_{x^*}\right\}$, there exists $u_0\neq0$ such that $f(z^*+u_0e_i)<f(z^*)$, which togehter with $z^*+u_0e_i\in F$ contradicts the optimality of $z^*$. Similar to the analysis in C\emph{ase 1}, we can obtain that $\nabla_if(z^*)=g'_i(0)=0$ for any $i\in\left(\Gamma^*_1\setminus\left\{\gamma_{x^*}\right\}\right)\cup I_2$. Thus, by Theorem \ref{thm-N-B} we deduce that $z^*$ is an ${\rm N^B}$-stationary point when $\Vert x^*\Vert_0=s$ and $\Vert y^*\Vert_0<t$.

\emph{Case 3}: $\Vert x^*\Vert_0<s$ and $\Vert y^*\Vert_0=t$.
For any $i\in\left\{\gamma_{x^*}+1,\ldots,m\right\}\cup\Gamma^*_2$, we have $z^*+ue_i\in F$ for any $u\in\mathbb{R}$. If $g'_i(0)\neq0$ for any $i\in\left\{\gamma_x+1,\ldots,m\right\}\cup\Gamma^*_2$, there exists $u_0\neq0$ such that $f(z^*+u_0e_i)<f(z^*)$, which together with $z^*+u_0e_i\in F$ contradicts the optimality of $z^*$. Similar to \emph{Case 1}, we can obtain that $\nabla_if(z^*)=g'_i(0)=0$ for any $i\in\left\{\gamma_{x^*}+1,\ldots,m\right\}\cup\Gamma^*_2$. Thus, by Theorem \ref{thm-N-B} we deduce that $z^*$ is an ${\rm N^B}$-stationary point when $\Vert x^*\Vert_0<s$ and $\Vert y^*\Vert_0=t$.

\emph{Case 4}: $\Vert x^*\Vert_0<s$ and $\Vert y^*\Vert_0<t$. For any $i\in \left\{i\in I:i>\gamma_{x^*}\right\}$, we have $z^*+ue_i\in F$ for any $u\in\mathbb{R}$. If $g'_i(0)\neq0$ for any $i\in \left\{i\in I:i>\gamma_{x^*}\right\}$, there exists $u_0\neq0$ such that $f(z^*+u_0e_i)<f(z^*)$, which together with $z^*+u_0e_i\in F$ contradicts the optimality of $z^*$. Similar to {\bf Case 1}, we can obtain that $\nabla_if(z^*)=g'_i(0)=0$ for any $i\in \left\{i\in I:i>\gamma_{x^*}\right\}$. Thus, by Theorem \ref{thm-N-B} we conclude that $z^*$ is an ${\rm N^B}$-stationary point when $\Vert x^*\Vert_0<s$ and $\Vert y^*\Vert_0<t$.

Combining \emph{Case 1} with \emph{Case 2}, we complete the proof of this theorem.
\ep

\subsubsection{First-order optimality conditions  by  Clarke tangent and normal cones}

Firstly, we need to find out the precise form of Clarke tangent and normal cones of $F$ at a point in $F$. The following two lemmas can be used to achieve our goal.

\begin{lemma}\emph{\cite[Theorem 2.4.5]{C1990} \label{NT-times}}
Let $Z=X\times Y$ where $X,Y$ are Banach spaces, and let $z=(x,y)\in C_1\times C_2$, where $C_1\subseteq X$, $C_2\subseteq Y$. Then
\begin{eqnarray*}
T_{C_1\times C_2}(z)=T_{C_1}(x)\times T_{C_2}(y)\quad\mbox{\rm and}\quad
N_{C_1\times C_2}(z)=N_{C_1}(x)\times N_{C_2}(y).
\end{eqnarray*}
\end{lemma}

\begin{lemma} \emph{\cite[Theorem 2.2]{PXZ2015} \label{CNT-cone-S2}}
For any $y\in F_2$ with $\Gamma_2:={\rm supp}(y)$, the Clarke tangent cone and corresponding normal cone of $F_2$ at $y$ are given respectively as
\begin{eqnarray}
T^{C}_{F_2}(y)&=&\left\{d\in\mathbb{R}^n:\ {\rm supp}(d)\subseteq\Gamma_2\right\}={\rm span}\left\{e_i,\ i\in\Gamma_2\right\},\label{ctc-y}\\
N^{C}_{F_2}(y)&=&\left\{d\in\mathbb{R}^n:\ d_i=0,\ \forall i\in\Gamma_2\right\}={\rm span}\left\{e_i,\ i\notin\Gamma_2\right\}.  \label{cnc-y}
\end{eqnarray}
\end{lemma}

To express the Clarke tangent and normal cones of $F$ at a point in $F$, by Lemmas \ref{NT-times} and \ref{CNT-cone-S2}, it is sufficient to represent the Clarke tangent  and normal cones of $F_1$ at a point in $F_1$, which are given as follows.

\begin{theorem}\label{CNT-cone-F1}
Let $F_1$ be defined by (\ref{fsbl}) and $\Gamma_1$ be defined by (\ref{supphuang}). Then the Clarke tangent cone and corresponding normal cone of $F_1$ at $x\in F_1$ can be written respectively as
\begin{eqnarray}
T^C_{F_1}(x)&=&\left\{d\in\mathbb{R}^{m}:\ {\rm supp}(d)\subseteq\Gamma_1\setminus\left\{\gamma_x\right\}\right\},\label{ctc}\\
N^C_{F_1}(x)&=&\left\{d\in\mathbb{R}^m:\ d_i=0,\ \forall i\in\Gamma_1\setminus\left\{\gamma_x\right\}\right\}.    \label{cnc}
\end{eqnarray}
\end{theorem}

\noindent{\bf Proof}. If (\ref{ctc}) holds, then by the definition of Clarke normal cone (see Definition \ref{Ccone}), it is easy to see that (\ref{cnc}) holds. Thus, it is enough to show that (\ref{ctc}) holds.
 Firstly, we prove that \begin{equation} \label{aabbcc} T^C_{F_1}(x)\subseteq\left\{d\in\mathbb{R}^{m}:\ {\rm supp}(d)\subseteq\Gamma_1\setminus\left\{\gamma_x\right\}\right\}.
  \end{equation} For any $d\in T^C_{F_1}(x)$, we will show $d\in \left\{d\in\mathbb{R}^{m}:\ {\rm supp}(d)\subseteq\Gamma_1\setminus\left\{\gamma_x\right\}\right\}$. For this purpose, we only need to prove ${\rm supp}(d)\subseteq\Gamma_1$ and $d_{\gamma_x}=0$. Since $d\in T^C_{F_1}(x)$, it follows from the definition of $T^{C}_{F_1}(x)$ that for any $\left\{x^k\right\}\subset F_1$ with ${\rm lim}_{k\to\infty}x^k=x$ and $\left\{\lambda_k\right\}\subset\mathbb{R}_+$ with ${\rm lim}_{k\to\infty}\lambda_k=0$, there exists $\left\{d^k\right\}$ with ${\rm lim}_{k\to\infty}d^k=d$ such that $x^k+\lambda_kd^k\in F_1$ for any $k\in\mathbb{N}$.
\begin{itemize}
\item[(i)] We prove ${\rm supp}(d)\subseteq\Gamma_1$. Suppose that ${\rm supp}(d)\not\subseteq\Gamma_1. $ Then there exists $i_0\in{\rm supp}(d)$ such that $i_0\notin\Gamma_1$, i.e., $x_{i_0}=0$ and $d_{i_0}\neq0$. Since ${\rm lim}_{k\rightarrow\infty}d^k=d$, one must have $d^k_{i_0}\rightarrow d_{i_0}$. Then there exists $k_0\in\mathbb{N}$ such that $d_{i_0}\neq0$ for all $k\ge k_0$. Choose $\left\{x^k\right\}\subset F_1$ such that ${\rm lim}_{k\rightarrow\infty}x^k=x$ and ${\rm supp}(x^k)=\Gamma_1\cup\Gamma_1^k$ with $|\Gamma_1\cup\Gamma^k_1|=s$, where $\Gamma_1^k\cap\Gamma_1=\emptyset$ and $i_0\notin\Gamma_1^k$. Trivially, $\Gamma_1^k=\emptyset$ if $|\Gamma_1|=s$. Let $\lambda_k$ be selected as
    \begin{eqnarray*}
    \lambda_k=\min_{i\in\Gamma_1\cup\Gamma_1^k}\frac{|x^k_i|}{|d^k_i|+k},
    \end{eqnarray*}
    then one has $\lambda_k>0$ for any $k\in\mathbb{N}$ and $\lambda_k\rightarrow 0$. For any $i\in\Gamma_1\cup\Gamma^k_1$, we have
    \begin{eqnarray*}
    |x^k_i+\lambda_kd^k_i|\ge|x^k_i|-\lambda_k|d^k_i|>0.
    \end{eqnarray*}
    Since $\lambda_k>0$, $x^k_{i_0}=0$ and $d^k_{i_0}\neq0$ for $k\ge k_0$, one must have $x^k_{i_0}+\lambda_kd^k_{i_0}\neq0$ for $k\ge k_0$. This, together $i_0\notin\Gamma_1\cup\Gamma_1^k$ and $|\Gamma_1\cup\Gamma^k_1|=s$, implies that $\Vert x^k+\lambda_kd^k\Vert_0>s$ for $k\ge k_0$. This contradicts $x^k+\lambda_kd^k\in F_1$. Thus ${\rm supp}(d)\subseteq\Gamma_1$.
\item[(ii)] We prove $d_{\gamma_x}=0$. Suppose that $d_{\gamma_x}\neq0$. It follows that for any $\left\{d^k\right\}$ with ${\rm lim}_{k\to\infty}d^k=d$, there exists $k_1\in\mathbb{N}$ such that $d^{k}_{\gamma_x}\neq0$ for any $k>k_1$. Set
  $
    \lambda_k=\frac{1}{|d^k_{\gamma_x}|+k} $  for all $ k>k_1,
   $
    then one has $\lambda_k>0$ for all $k>k_1$ and ${\rm lim}_{k\to\infty}\lambda_k=0$. For any $j<\gamma_x$, one must have $j\notin\Gamma_1$. It follows that $d_j=0$ by ${\rm supp}(d)\subseteq\Gamma_1$. Thus, ${\rm lim}_{k\to\infty}d^k_j=0$. Then for $\epsilon<1$, there exists $k_2\in\mathbb{N}$ such that $\lambda_kd^k_j<\epsilon$ for any $k>k_1$ and $j<\gamma_x$. By Lemma \ref{x-k-1}, there exists $k_0\in\mathbb{N}$ such that $x^k_j=0$ for any $k>k_0$ and $j<\gamma_x$ and $x^k_{\gamma_x}=1$. Let $K=\max\left\{k_0,k_1,k_2\right\}$, for any $k>K,j<\gamma_x$, one has $(x^k+\lambda_kd^k)_j<\epsilon$,
    \begin{eqnarray*}
    (x^k+\lambda_kd^k)_{\gamma_x}=1+\frac{d^k_{\gamma_x}}{|d^k_{\gamma_x}|+k}\neq0\ \text{and}\ (x^k+\lambda_kd^k)_{\gamma_x}=1+\frac{d^k_{\gamma_x}}{|d^k_{\gamma_x}|+k}\neq1,
    \end{eqnarray*}
    which contradicts $x^k+\lambda_kd^k\in F_1$. So $d_{\gamma_x}=0.$
\end{itemize}
Thus (\ref{aabbcc}) holds. We now  prove that \begin{equation} \label{ccbbaa} \left\{d\in\mathbb{R}^{m}:\ {\rm supp}(d)\subseteq\Gamma_1\setminus\left\{\gamma_x\right\}\right\}\subseteq T^C_{F_1}(x).
\end{equation}  To this goal, for any $d\in \left\{d\in\mathbb{R}^{m}:\ {\rm supp}(d)\subseteq\Gamma_1\setminus\left\{\gamma_x\right\}\right\}$, we need to show $d\in T^C_{F_1}(x)$.
For any $\left\{x^k\right\}\subset F_1$ such that ${\rm lim}_{k\rightarrow\infty}x^k=x$ and any $\left\{\lambda_k\right\}\subset\mathbb{R}_+$ such that ${\rm lim}_{k\to\infty}\lambda_k=0$, it follows that there exists $k_0\in\mathbb{N}$ such that $\Gamma_1\subseteq{\rm supp}(x^k)$ and $x^k_{\gamma_x}=1$ for any $k\ge k_0$. Moreover, we choose a sequence $\{d^k\}$  defined as
\begin{equation}\label{dk-C}
d^k:=\left\{
\begin{array}{ll}
0\ &\ \mbox{\rm if}\; k=1,2,\ldots,k_0,\\
x^k-x+d\ &\ \mbox{\rm if}\;k=k_0+1,k_0+2,\ldots,
\end{array}
\right.
\end{equation}
which satisfies that  ${\rm lim}_{k\to\infty}d^k=d$ by ${\rm lim}_{k\to\infty}x^k=x$. For $k<k_0$, we have
\begin{align*}
d^k_j &  =  x_j^k-x_j+d_j=0-0+0,\ \forall j<\gamma_x,\\
d^k_{\gamma_x} & =x^k_{\gamma_x}-x_{\gamma_x}+d_{\gamma_x}=1-1+0=0,\\
{\rm supp}(d^k) & \subseteq{\rm supp}(x^k)\cup{\rm supp}(x)\cup{\rm supp}(d)={\rm supp}(x^k).
\end{align*}
So we have $\Vert x^k+\lambda_kd^k\Vert_0\le s$ and $x^k_{\gamma_x}+\lambda_kd^k_{\gamma_x}=1$ for any $k>k_0$. These, together with $x^k+\lambda_kd^k=x^k\in F_1$ for any $k\le k_0$, imply that $x^k+\lambda_kd^k\in F_1$ for any $k\in\mathbb{N}$. Therefore, $d\in T^C_{F_1}(x)$, and hence (\ref{ccbbaa}) holds. Thus (\ref{ctc}) holds. The proof  is complete.
\ep

Now, we provide the precise forms of the Clarke tangent cone and Clarke normal cone of $F$ at a point in $F$.
\begin{theorem}
For any $z=(x,y)\in F$ with $F$ being defined by (\ref{sprbl-set}),
the Clarke tangent cone and corresponding normal cone of $F$ at $z$ can be written as
\begin{eqnarray}
T^C_{F}(z)&=&\left\{d\in\mathbb{R}^{m+n}:\ {\rm supp}(d)\subseteq\Gamma\setminus\left\{\gamma_x\right\}\right\}, \label{ctc-z}\\
N^C_{F}(z)&=&\left\{d\in\mathbb{R}^{m+n}:\ d_i=0,\ \forall i\in\Gamma\setminus\left\{\gamma_x\right\}\right\},\label{cnc-z}
\end{eqnarray}
respectively, where $\Gamma$ is defined by (\ref{supphuang}).
\end{theorem}

This result follows immediately from Theorem \ref{CNT-cone-F1} and Lemmas \ref{NT-times} and \ref{CNT-cone-S2}. Next, we provide a sufficient and necessary condition for  ${\rm N^C}$-stationary  and ${\rm T^C}$-stationary points.

\begin{theorem}\label{thm-N-C}
Let $z^*=(x^*,y^*)\in F$, where $F$ is defined by (\ref{sprbl-set}), and let  $\Gamma^*$ be defined by (\ref{gm-star}). Then (i) $z^*$ is an ${\rm N^C}$-stationary point if and only if $\nabla f(z^*)\in \mathbb{R}^{m+n}$ satisfies
\begin{equation}\label{N-C-P}
\nabla_if(z^*)=0,\;\forall i\in\Gamma^*\setminus\left\{\gamma_{x^*}\right\}.
\end{equation}
(ii)
$z^*$ is a ${\rm T^C}$-stationary point if and only if (\ref{N-C-P}) holds.
\end{theorem}

\noindent{\bf Proof}. Item (i) follows directly from (\ref{cnc-z}). We now show Item (ii).  By Definition \ref{def-T-stp}, we need to show that $\nabla^C_Ff(z^*)=0$ if and only if (\ref{N-C-P}) holds.

($\Rightarrow$)  Suppose there exists $i_0\in\Gamma\setminus\left\{\gamma_{x^*}\right\}$ such that $\nabla_{i_0}f(z^*)\neq 0$. Pick $\bar{d}=(\bar{d}_1,\bar{d}_2)\in\mathbb{R}^{m+n}$ such that $\bar{d}_{i_0}=-\nabla_{i_0}f(z^*)$ and other entries are zero. Since ${\rm supp}(\bar{d})\subseteq \Gamma^*\setminus\left\{\gamma_{x^*}\right\}$, it follows from (\ref{ctc-z}) that $\bar{d}\in T^C_F(z^*)$. Then, by Definition \ref{def-T-stp}, we have
$$
\Vert\nabla_F^Cf(z^*)+\nabla f(z^*)\Vert^2\le\Vert\bar{d}+\nabla f(z^*)\Vert^2.
$$
Since $\nabla^C_Ff(z^*)=0$, the above inequality is equivalent to (\ref{i0le}), which further implies $(\nabla_{i_0}f(z^*))^2\le0$. This contradicts the assumption $\nabla_{i_0}f(z^*)\neq 0$.

($\Leftarrow$) Suppose there exists $\bar{d}\in T^C_F(z^*)$ with $\bar{d}\neq0$ such that (\ref{bar-d-le}) holds.
It follows from (\ref{ctc-z}) that ${\rm supp}(\bar{d})\subseteq\left(\Gamma^*\setminus\left\{\gamma_{x^*}\right\}\right)$. Thus (\ref{bar-d-le}) is equivalent to
$$
\sum_{i\in\Gamma^*\setminus\left\{\gamma_{x^*}\right\}}(\bar{d}_i+\nabla_if(z^*))^2\le\sum_{i\in\Gamma^*\setminus\left\{\gamma_{x^*}\right\}}(\nabla_if(z^*))^2.
$$
This, together with (\ref{N-C-P}), implies that $\sum_{i\in\Gamma^*\setminus\left\{\gamma_{x^*}\right\}}(\bar{d}_i)^2\le0$, a contradiction.
\ep

By Theorems \ref{thm-N-B} and \ref{thm-T-B}, we immediately obtain the following Corollary.

\begin{corollary}
Let $z^*=(x^*,y^*)\in F,$ where $F$ is defined by (\ref{sprbl-set}). Then
$z^*$ is an ${\rm N^C}$-stationary point of (\ref{originalmodel}) if and only if it is a ${\rm T^C}$-stationary point of (\ref{originalmodel}).
\end{corollary}

We now discuss the relation between ${\rm N^B}$-stationary point and ${\rm N^C}$-stationary point.
\begin{theorem}\label{B-to-C}
Let $z^*=(x^*,y^*)\in F,$ where $F$ is defined by (\ref{sprbl-set}).
If $z^*$ is an ${\rm N^B}$-stationary point of (\ref{originalmodel}), then $z^*$ is also an ${\rm N^C}$-stationary point.
\end{theorem}

\noindent{\bf Proof}.  Let $\Gamma^*,\Gamma^*_1,\Gamma^*_2$ be defined as in (\ref{gm-star}). It is sufficient to consider the following four cases:
\begin{itemize}
\item[(i)]  $\Vert x^*\Vert_0=s$ and $\Vert y^*\Vert_0=t$. For this case, the result can be obtained directly from (\ref{N-B-c1}) and (\ref{N-C-P}).

\item[(ii)]   $\Vert x^*\Vert_0<s$ and $\Vert y^*\Vert_0=t$. Since $\Gamma^*_1\setminus\left\{\gamma_{x^*}\right\}\subseteq\left\{\gamma_{x^*}+1,\ldots,m\right\}$, it follows that $\nabla_if(z^*)=0$ for any $i\in\Gamma^*_1\setminus\left\{\gamma_{x^*}\right\}$ if $\nabla_if(z^*)=0$ for any $i\in\left\{\gamma_{x^*}+1,\ldots,m\right\}$. Thus it follows from (\ref{bnc-z}) and (\ref{cnc-z}) that $z^*$ is an ${\rm N^C}$-stationary point if $z^*$ is an ${\rm N^B}$-stationary point.

\item[(iii)]   $\Vert x^*\Vert_0=s$ and $\Vert y^*\Vert_0<t$. Since $\Gamma^*_1\setminus\left\{\gamma_{x^*}\right\}\subseteq(\Gamma^*_1\cup I_2)\setminus\left\{\gamma_{x^*}\right\}$, similar to the analysis of Case (ii), it follows from (\ref{bnc-z}) and (\ref{cnc-z}) that $z^*$ is an ${\rm N^C}$-stationary point if $z^*$ is an ${\rm N^B}$-stationary point.

\item[(iv)]  $\Vert x^*\Vert_0<s$ and $\Vert y^*\Vert_0<t$. Since $\Gamma^*_1\setminus\left\{\gamma_{x^*}\right\}\subseteq\left\{i\in I:\ i>\gamma_{x^*}\right\}$, it follows from (\ref{bnc-z}) and (\ref{cnc-z}) that $z^*$ is an ${\rm N^C}$-stationary point if $z^*$ is an ${\rm N^B}$-stationary point.
\end{itemize}
Thus an N$^B$-stationary point must be an N$^C$-stationary point.
\ep

\begin{corollary}
If $|\Gamma^*_1|=s$ and $|\Gamma^*_2|=t$, then the $N^B$-stationary point and $T^B$-stationary point are equivalent to the $N^C$-stationary point and $T^C$-stationary point, respectively.
\end{corollary}

Finally, we discuss the relation between ${\rm N^C}$-stationary point (${\rm T^C}$-stationary point) of (\ref{originalmodel}) and the optimal solution of (\ref{originalmodel}).
 The following result is implied directly from Theorems \ref{B-nec} and  \ref{B-to-C}.

\begin{theorem}\label{N-C-optim}
$z^*=(x^*,y^*)$ is an optimal solution of (\ref{originalmodel}), then $z^*$ is an ${\rm N^C}$-stationary point, and $z^*$ is also a ${\rm T^C}$-stationary point.
\end{theorem}

\subsection{Coordinate-Wise minimum}

Motivated by the concept of coordinate-wise minimum (CW minimum) for the problem (\ref{nonlinear-model}) proposed in \cite{BE2013}, we extend this notion to the problem (\ref{originalmodel}). Because of the equality constraint $e^{\top}_{\gamma_x}x=1$,  the definition is appropriately modified to fit the current framework.

\begin{definition}\label{cw-min}
Let $F$ be defined by (\ref{sprbl-set}), then $z^*$ is called a CW minimum of (\ref{originalmodel}) if and only if one of the following conditions is satisfied:
\begin{itemize}
\item[(i)] When $\Vert x^*\Vert_0=s$ and $\Vert y^*\Vert_0=t$, for any index pair $(i,j)\in I_1\times I_1$ or $I_2\times I_2$ with $i\in\Gamma^*\setminus\left\{\gamma_{x^*}\right\}$ and $j>\gamma_x$, one has
\begin{eqnarray}\label{cw-eq2}
f(z^*)=\min_{u\in\mathbb{R}}f(z^*-z^*_ie_i+ue_j).
\end{eqnarray}
\item[(ii)] When $\Vert x^*\Vert_0<s$ and $\Vert y^*\Vert_0<t$, for any $i>\gamma_{x^*}$, one has
\begin{eqnarray}\label{cw-eq1}
f(z^*)=\min_{u\in\mathbb{R}}f(z^*+ue_i).
\end{eqnarray}
\item[(iii)] When $\Vert x^*\Vert_0=s$ and $\Vert y^*\Vert_0<t$, for any $i\in I_2$, one has (\ref{cw-eq1}) holds, and for any $(i,j)\in I_1\times I_1$ with $i\in\Gamma^*_1\setminus\left\{\gamma_{x^*}\right\}$ and $j>\gamma_{x^*}$, one has (\ref{cw-eq2}) holds.
\item[(iv)] When $\Vert x^*\Vert_0<s$ and $\Vert y*\Vert_0=t$, for any $i\in I_1$ with $i>\gamma_{x^*}$, one has (\ref{cw-eq1}) holds, and for any $i\in\Gamma^*_2$ and $j\in I_2$, one has (\ref{cw-eq2}) holds.
\end{itemize}
\end{definition}

Based on the above definition, it follows that any optimal solution of the problem must be a CW minimum, as stated in the following theorem.

\begin{theorem}\label{cw-nec}
Let $z^*$ be an optimal solution of (\ref{originalmodel}), then $z^*$ is a CW minimum of (\ref{originalmodel}).
\end{theorem}
\noindent{\bf Proof.} It follows directly from the optimality of $z^*$.
\ep

It can be observed that a CW minimum is also an $N^B$-stationary point, which is stated in the following theorem.

\begin{theorem}
Let $z^*=(x^*,y^*)\in F,$ where $F$ is defined by (\ref{sprbl-set}). If $z^*$ is a CW minimum of (\ref{originalmodel}), then $\nabla f(z^*)\in \mathbb{R}^{m+n}$ satisfies (\ref{N-B-c1}), and hence, $z^*$ is an $N^B$-stationary point.
\end{theorem}

\noindent{\bf Proof.} Assume that $z^*\in F$ is a CW minimum of (\ref{originalmodel}). We first show that $\nabla f(z^*)\in \mathbb{R}^{m+n}$ satisfies (\ref{N-B-c1}), which is divided into the following four cases.

(i) $\Vert x^*\Vert_0=s$ and $\Vert y^*\Vert_0=t$. By Definition \ref{cw-min}, it follows that (\ref{cw-eq2}) is satisfied for any index pair $(i,j)\in I_1\times I_1$ or $I_2\times I_2$ with $i\in\Gamma^*\setminus\left\{\gamma_{x^*}\right\}$ and $j>\gamma_x$,
then we have
$$
f(z^*)=\min_{u\in\mathbb{R}}f(z^*-z^*_ie_i+ue_i),
$$
which becomes
\begin{eqnarray}\label{v-optim}
f(z^*)=\min_{v\in\mathbb{R}}f(z^*+ve_i),
\end{eqnarray}
by denoting $v:=u-z^*_i$. It follows from (\ref{v-optim}) that $v=0$ is an optimal solution of the problem $\min_{v\in\mathbb{R}}h_i(v)=f(z^*+ve_i)$. By the optimality condition $h_i'(0)=0$ for any $i\in\Gamma^*\setminus\left\{\gamma_{x^*}\right\}$, we conclude that for any CW minimum $z^*$ of (\ref{originalmodel}) with $\Vert x^*\Vert_0=s$ and $\Vert y^*\Vert_0=t$, we have
$$
\nabla_i f(z^*)=0,\ \forall i\in\Gamma^*\setminus\left\{\gamma_{x^*}\right\},
$$
which implies that case $\Vert x^*\Vert_0=s$ and $\Vert y^*\Vert_0=t$ in (\ref{N-B-c1}) is satisfied. That is, $z^*$ is an $N^B$-stationary point in this case.

(ii) $\Vert x^*\Vert_0<s$ and $\Vert y^*\Vert_0<t$. It follows from Definition \ref{cw-min} that (\ref{v-optim}) is already satisfied for any $i>\gamma_{x^*}$. By the same analysis in (i), we can obtain that
$$
\nabla_i f(z^*)=0,\ \forall i>\gamma_{x^*},
$$
which implies that case $\Vert x^*\Vert_0<s$ and $\Vert y^*\Vert_0<t$ in (\ref{N-B-c1}) is satisfied. That is, $z^*$ is an $N^B$-stationary point in this case.

(iii) $\Vert x^*\Vert_0=s$ and $\Vert y^*\Vert_0<t$. By Definition \ref{cw-min}, one has (\ref{cw-eq2}) holds for any $(i,j)\in I_1\times I_1$ with $i\in\Gamma^*_1\setminus\left\{\gamma_{x^*}\right\}$ and $j>\gamma_{x^*}$. By similar analysis in (i), we can obtain that
\begin{eqnarray}\label{cw-eq3}
\nabla_i f(z^*)=0,\ \forall i\in\Gamma^*_1\setminus\left\{\gamma_{x^*}\right\}.
\end{eqnarray}
By Definition \ref{cw-min}, one has (\ref{cw-eq1}) holds for any $i\in I_2$. Then by similar analysis in (ii), we have $\nabla_if(z^*)=0$ for any $i\in I_2$, which implies that case $\Vert x^*\Vert_0=s$ and $\Vert y^*\Vert_0<t$ in (\ref{N-B-c1}) is satisfied. That is, $z^*$ is an $N^B$-stationary point in this case.

(iv) $\Vert x^*\Vert_0<s$ and $\Vert y^*\Vert_0=t$. Since for any $i\in I_1$ with $i>\gamma_{x^*}$, one has (\ref{cw-eq1}) holds, and for any $i\in\Gamma^*_2$ and $j\in I_2$, one has (\ref{cw-eq2}) holds. Analogous to (iii), we can conclude that $\nabla_if(z^*)=0$ for any $i\in\left\{\gamma_{x^*}+1,\ldots,m\right\}\cup\Gamma^*_2$, which implies that case $\Vert x^*\Vert_0<s$ and $\Vert y^*\Vert_0=t$ in (\ref{N-B-c1}) is satisfied. That is, $z^*$ is an $N^B$-stationary point in this case.

Combining the above cases (i)-(iv), we obtain that $\nabla f(z^*)\in \mathbb{R}^{m+n}$ satisfies (\ref{N-B-c1}). Furthermore, by Theorem \ref{thm-N-B} we obtain that $z^*$ is an $N^B$-stationary point.
\ep

Note that an $N^B$-stationay point may not be a CW minimum, we present the following example to illustrate it.

\begin{example}
Consider the problem (\ref{originalmodel}) with $b=(5,0)^\top$, tensor $\mathscr{A}=(a_{j_1j_2j_3})\in\mathbb{R}^{2\times3\times3}$, where $a_{111}=1, a_{112}=1, a_{113}=3, a_{121}=1, a_{122}=1, a_{133}=1, a_{211}=1, a_{213}=-1, a_{222}=1, a_{231}=1$ and all other entries are zero. Then for any $z=(x,y)\in\mathbb{R}^3\times \mathbb{R}^3$, the objective function in (\ref{originalmodel}) is defined by $f(z):=\frac{1}{2}\Vert \mathscr{A}xy-b\Vert^2$ with
\begin{equation*}
\mathscr{A}xy-b=\left[
\begin{array}{l}
x_1y_1+x_1y_2+3x_1y_3+x_2y_1+x_2y_2+x_3y_3-5\\
x_1y_1-x_1y_3+x_2y_2+x_3y_1
\end{array}
\right].
\end{equation*}
Let $s=t=2$ and $\bar{z}=(1,1,0,1,1,0)^\top$. We have $\Vert\bar{x}\Vert_0=s$, $\Vert\bar{y}\Vert_0=t$, $\gamma_{\bar{x}}=1$, $\bar{\Gamma}=\{1,2,4,5\}$, $f(\bar{z})=\frac{5}{2}$ and $\nabla f(\bar{z})=(0,0,2,0,0,-5)^\top$. It is easy to see that $\nabla_if(\bar{z})=0$ for $i\in \bar{\Gamma}\setminus \{\gamma_{\bar{x}}\}$. Thus, by Theorem \ref{thm-N-B} we obtain that that $\bar{z}$ is an $N^B$-stationary point. However, if we define $z^*:=\bar{z}-e_5+e_6=(1,1,0,1,0,1)^\top$, then we have $f(z^*)=0$. Obviously, $f(\bar{z})>f(z^*)$, which implies that (\ref{cw-eq2}) is not satisfied at $\bar{z}$. In another word, $\bar{z}$ is not a CW minimum.
\end{example}

\subsection{$L$-like stationary point}

Recall that for the single-block variable sparse optimization problem (\ref{nonlinear-model}), a vector $x^*\in C_s:=\left\{x\in\mathbb{R}^m:\  \Vert x\Vert_0 \le s\right\}$ is said to be an  $L$-stationary point if and only if $x^*\in{\bf P}_{C_s}(x^*-\nabla f(x^*)/L)$ for some $L>0$, where ${\bf P}_{C_s}(x):={\rm argmin}_{u\in C_s}\Vert u-x\Vert$ (see, e.g., \cite{BE2013}). Since the problem (\ref{originalmodel}) has a special constraint $e_{\gamma_x}^\top x=1$, it is difficult for us to characterize the $L$-stationary point for this problem by using the projection operator ${\bf P}_F$, where $F$ is defined by (\ref{sprbl-set}). In this section, similar to the idea used in the $L$-stationary point, we introduce and discuss the concept of an {\it $L$-like stationary point}. For this purpose, we introduce the following concept of {\it like-projection operator}.
\par
Before introducing this new concept, let us first give the definition below.
\begin{definition}\label{psiz}
Let $\psi_z:\mathbb{R}^{m+n}\to\mathbb{R}^{m+n}$ be a mapping defined as follows: given $z=(x,y)\in\mathbb{R}^{m+n}$, for any $u\in\mathbb{R}^{m+n}$,
\begin{itemize}
\item[(a)] if $x=0$, $\psi_z(u):=u$;
\item[(b)] if $x\neq0$,
\begin{equation*}
\left(\psi_z(u)\right)_i:=\left\{
\begin{array}{ll}
x_{\gamma_u^z}u_i&\ \mbox{\rm if}\; i\in I_1,\\
\frac{u_i}{x_{\gamma_u^z}}&\ \mbox{\rm if}\; i\in I_2,
\end{array}
\right.
\end{equation*}
where $\gamma_u^z$ is given by
\begin{eqnarray}\label{gammamin}
\gamma_u^z:=\min\{i\in I_1:u_i\neq0, z_{i}\neq 0\}.
\end{eqnarray}
\end{itemize}
\end{definition}
By using the mapping $\psi_z(u)$ above, we can now define the like-projection.
\begin{definition}\label{def-proj}
For any $z=(x,y)\in\mathbb{R}^{m+n}$, let $F$ be defined by (\ref{sprbl-set}). The like-projection of $z$ onto $F$ is defined by
$$
{\bar {\bf P}}_F(z)=\arg\min_{u\in F}\Vert\psi_z(u)-z\Vert^2.
$$
\end{definition}

\begin{remark}\label{remark3.1}
The like-projection operator ${\bar {\bf P}}_F(\cdot)$ introduced in Definition \ref{def-proj} is determined by the set $F$ defined by (\ref{sprbl-set}), it is different from the classical projection operator ${\bf P}_F(\cdot)$ (i.e., ${\bf P}_F(z):={\rm argmin}_{u\in F}\Vert u-z\Vert$ for $z\in \mathbb{R}^{m+n}$), but when $z_{\gamma_u^z}=1$ or $x_i=0$ for all $i\in I_1$, the like-projection of $z$ onto $F$ coincides with the projection of $z$ onto $F$.
\end{remark}

It is easy to find the projection of an arbitrary point $x\in \mathbb{R}^m$ onto the sparse set $C_s=\left\{x\in\mathbb{R}^m:\ \Vert x\Vert_0 \le s\right\}$. A natural question is how to calculate the like-projection of a point $z\in \mathbb{R}^{m+n}$ onto the set $F$ defined by (\ref{sprbl-set}). The following lemma addresses this question.

\begin{lemma}\label{lem3-4}
For any $z=(x,y)\in\mathbb{R}^{m+n}$,  let $q_1,q_2$ be the smallest integers and $r_1,r_2$ be the largest integers such that
$$
 M_{q_1}(|x|)=M_{q_1+1}(|x|)=\cdots=M_{s}(|x|)=\cdots=M_{r_1}(|x|),$$
 $$M_{q_2}(|y|)=M_{q_2+1}(|y|)=\cdots=M_{t}(|y|)=\cdots=M_{r_2}(|y|),
$$
and  $\Lambda_1 $ and $\Lambda_2 $  are defined as
\begin{eqnarray*}
\Lambda_1:=\{i\in I_1: |x_i|=M_{s}(|x|)\},\;\;
\Lambda_2:=\{i\in I_2: |z_i|=M_{t}(|y|)\}.
\end{eqnarray*}
If $u^*\in {\bar {\bf P}}_F(z)$, then there exist sets $\Upsilon_1\subseteq\Lambda_1$ and $\Upsilon_2\subseteq\Lambda_2$ with $|\Upsilon_1|=s-q_1+1$ and $|\Upsilon_2|=t-q_2+1$ such that one of the following holds:
\begin{itemize}
\item $x\neq 0$ and $u^*\in\mathbb{R}^{m+n}$ is defined by
\begin{equation}\label{P-fml}
u_i^*=\left\{
\begin{array}{ll}
\frac{z_i}{x_{\gamma_{\rm min}^*}}&\ \mbox{\rm if}\; i\in \Omega_1:=\{i\in I_1: |z_i|>M_s(|x|)\}\bigcup\Upsilon_1,\\
x_{\gamma_{\rm min}^*}z_i&\ \mbox{\rm if}\; i\in \Omega_2:=\{i\in I_2: |z_i|>M_t(|y|)\}\bigcup\Upsilon_2,\\
0&\ \mbox{\rm otherwise},
\end{array}
\right.
\end{equation}
where $\gamma_{\rm min}^*=\min\{i\in\Omega_1:x_i\neq 0\}$.
\item $x=0$ and $u^*\in\mathbb{R}^{m+n}$ is defined by
\begin{eqnarray}\label{P-fmle1}
\begin{array}{rcl}
u_{i_0}^*&=&1\; \mbox{\rm for some}\ i_0\in \Upsilon_1=I_1,\\
u_i^*&=&\left\{
\begin{array}{ll}
0&\ \mbox{\rm if}\; i\in I_1\setminus\{i_0\},\\
z_i&\ \mbox{\rm if}\; i\in \Omega_2:=\{i\in I_2: |z_i|>M_t(|y|)\}\bigcup\Upsilon_2,\\
0&\ \mbox{\rm if}\; i\in I_2\setminus\Omega_2.
\end{array}\right.
\end{array}
\end{eqnarray}
\end{itemize}
Conversely, for any $\Upsilon_1\subseteq\Lambda_1$ with $|\Upsilon_1|=s-q_1+1$ and $\Upsilon_2\subseteq\Lambda_2$ with $|\Upsilon_2|=t-q_2+1$, if $u^*$ is given by (\ref{P-fml}) for $x\neq0$ or $u^*$ is given by  (\ref{P-fmle1}) for $x=0$, then $u^*\in\bar{\bf P}_F(z)$.
\end{lemma}

\noindent{\bf Proof}.
($\Rightarrow$)  We first show that either (\ref{P-fml}) or (\ref{P-fmle1}) holds if $u^*\in {\bar {\bf P}}_F(z)$.  It suffices to consider the following two cases.

\emph{Case 1:} $z\in\mathbb{R}^{m+n}$ with $x\neq0$. In this case, we show that (\ref{P-fml}) holds if $u^*\in {\bar {\bf P}}_F(z)$. For any $u=(u^1,u^2)\in F$, it follows from the definition of $F$ that $\|u^1\|_0\leq s$ and $\|u^2\|_0\leq t$, and hence,
there exist sets $\Upsilon_1\subseteq\Lambda_1$ and $\Upsilon_2\subseteq\Lambda_2$ with $|\Upsilon_1|=s-q_1+1$ and $|\Upsilon_2|=t-q_2+1$ such that
$$
\Omega_1:=\{i\in I_1: |z_i|>M_s(|x|)\}\bigcup\Upsilon_1\;\; \mbox{\rm and}\;\; \Omega_2:=\{i\in I_2: |z_i|>M_t(|y|)\}\bigcup\Upsilon_2
$$
satisfy $\mbox{\rm supp}(u^1)=\mbox{\rm supp}(\{z_i: i\in \Omega_1\})$ and $\mbox{\rm supp}(u^2)=\mbox{\rm supp}(\{z_i: i\in \Omega_2\})$. It follows from Definition \ref{psiz} that  $\psi_z(u)=(\psi_z^1(u),\psi_z^2(u^2))$ with $\psi_z^1(u)=x_{\gamma_u^z}u^1$ and $\psi_z^2(u)=\frac{u^2}{x_{\gamma_u^z}}$, where $\gamma_u^z$ is defined as (\ref{gammamin}). Thus,
\begin{eqnarray}\label{phi-u-eq}
\Vert \psi_z(u)-z\Vert^2
&=&\sum_{i\in I_1}|\left(\psi_z^1(u)\right)_i-z_i|^2+\sum_{j\in I_2}|\left(\psi_z^2(u)\right)_j-z_j|^2\nonumber\\
&=&\sum_{i\in\Omega_1}|x_{\gamma_{u}^z}u^1_i-z_i|^2 +\sum_{j\in\Omega_2}|u^2_j/x_{\gamma_{u}^z}-z_j|^2\nonumber\\
&&+\sum_{i\in I_1\setminus\Omega_1}|z_i|^2+\sum_{j\in I_2\setminus\Omega_2}|z_j|^2.
\end{eqnarray}
Furthermore, we have
\begin{eqnarray}\label{ustar-z-e1}
\min_{u\in F}\Vert \psi_z(u)-z\Vert^2
=\min_{u\in F}\left\{\sum_{i\in\Omega_1}|x_{\gamma_{u}^z}u^1_i-z_i|^2 +\sum_{j\in\Omega_2}|u^2_j/x_{\gamma_{u}^z}-z_j|^2\right\}.
\end{eqnarray}
Since $u^*\in {\bar {\bf P}}_F(z)={\rm argmin}_{u\in F}\Vert\psi_z(u)-z\Vert^2$, it follows from (\ref{ustar-z-e1}) that
\begin{equation}\label{ustar-z-e2}
u_i^*=\left\{
\begin{array}{ll}
\frac{z_i}{x_{\gamma_{u^*}^z}}&\ \mbox{\rm if}\; i\in \Omega_1,\\
x_{\gamma_{u^*}^z}z_i&\ \mbox{\rm if}\; i\in \Omega_2,\\
0&\ \mbox{\rm otherwise}.
\end{array}
\right.
\end{equation}
From (\ref{ustar-z-e2}), it is easy to see that for any $i\in I$, $u_i^*\neq 0$ if $z_i\neq 0$, and hence, $\gamma_{\rm min}^*=\gamma_{u^*}^z$ where
$\gamma_{\rm min}^*=\min\{i\in\Omega_1:x_i\neq 0\}$. Thus, (\ref{ustar-z-e2}) is nothing but (\ref{P-fml}). So   (\ref{P-fml}) holds for this case.

\emph{Case 2}:  $z\in\mathbb{R}^{m+n}$ with $x=0$. In this case, we show that (\ref{P-fmle1}) holds when $u^*\in {\bar {\bf P}}_F(z)$. For any $u=(u^1,u^2)\in F$, it follows from the definition of $ F$ that there exist $i_0\in\Upsilon_1=I_1$ such that $u_{i_0}^1=1$ and a set $\Upsilon_2\subseteq\Lambda_2$ with $|\Upsilon_2|=t-q_2+1$ such that
$
\Omega_2:=\{i\in I_2: |z_i|>M_t(|y|)\}\cup\Upsilon_2
$
satisfies $\mbox{\rm supp}(u^2)=\mbox{\rm supp}(\{z_i: i\in \Omega_2\})$. It follows from Definition \ref{psiz} that  $\psi_z(u)=u$. Thus, in a manner similar to (\ref{phi-u-eq}), one has
\begin{eqnarray*}
\Vert \psi_z(u)-z\Vert^2
=1+\sum_{i\in I_1\setminus\{i_0\}}|u^1_i-z_i|^2 +\sum_{j\in\Omega_2}|u^2_j-z_j|^2+\sum_{j\in I_2\setminus\Omega_2}|z_j|^2.
\end{eqnarray*}
Furthermore, we have
\begin{eqnarray}\label{ustar-z-e3}
\min_{u\in F}\Vert \psi_z(u)-z\Vert^2
=\min_{u\in F}\left\{\sum_{i\in I_1\setminus\{i_0\}}|u^1_i-z_i|^2 +\sum_{j\in\Omega_2}|u^2_j-z_j|^2\right\}.
\end{eqnarray}
Since $u^*\in {\bar {\bf P}}_F(z)={\rm argmin}_{u\in F}\Vert\psi_z(u)-z\Vert^2$, it follows from (\ref{ustar-z-e3}) and $u_{i_0}^{1*}=1$ that
$
u_{i_0}^*=1 $ for some $ i_0\in \Upsilon_1=I_1,$ and
$$ u_i^*=\left\{
\begin{array}{ll}
0&\ \mbox{\rm if}\; i\in I_1\setminus\{i_0\},\\
z_i&\ \mbox{\rm if}\; i\in \Omega_2:=\{i\in I_2: |z_i|>M_t(|y|)\}\bigcup\Upsilon_2,\\
0&\ \mbox{\rm if}\; i\in I_2\setminus\Omega_2.
\end{array}\right.
$$
Thus (\ref{P-fmle1}) holds for this case.

($\Leftarrow$) We now show that $u^*\in {\bar {\bf P}}_F(z)$ when (\ref{P-fml}) or (\ref{P-fmle1}) holds. It is also sufficient to consider two cases.

\emph{Case 1}:  $z\in\mathbb{R}^{m+n}$ with $x\neq0$. In this case, let $u^*=(u^{1*},u^{2*})\in \mathbb{R}^m\times\mathbb{R}^n$ be defined by (\ref{P-fml}), then $\gamma^z_{u^*}=\min\left\{i\in\Omega_1:\ x_i\neq0\right\}=\gamma^*_{\rm min}$, where $\gamma_{u^*}^z$ is defined as (\ref{gammamin}). It follows from Definition \ref{psiz} that  $\psi_z(u^*)=(\psi_z^1(u^*),\psi_z^2(u^*))$ with $\psi_z^1(u^*)=x_{\gamma_{\min}^*}u^{1*}$ and $\psi_z^2(u^*)=\frac{u^{2*}}{x_{\gamma_{\min}}^*}$. Thus,
\begin{eqnarray}\label{ustar-z}
\Vert \psi_z(u^*)-z\Vert^2
&=&\sum_{i\in I_1}|u^*_ix_{\gamma_{\min}^*}-z_i|^2+\sum_{j\in I_2}|u^*_j/x_{\gamma_{\min}^*}-z_j|^2\nonumber\\
&=&\sum_{i\in\Omega_1}|u^*_ix_{\gamma_{\min}^*}-z_i|^2 +\sum_{j\in\Omega_2}|u^*_j/x_{\gamma_{\min}^*}-z_j|^2+\sum_{i\in I_1\setminus\Omega_1}|z_i|^2+\sum_{j\in I_2\setminus\Omega_2}|z_j|^2\nonumber\\
&=&\sum_{i\in\Omega_1}|z_i-z_i|^2+\sum_{j\in\Omega_2}|z_j-z_j|^2+\sum_{i\in I_1\setminus\Omega_1}|z_i|^2+\sum_{j\in I_2\setminus\Omega_2}|z_j|^2\nonumber\\
&=&\sum_{i\in I_1\setminus\Omega_1}|z_i|^2+\sum_{j\in I_2\setminus\Omega_2}|z_j|^2\nonumber\\
&=&\sum_{i=1}^{m-s}\left(M_{s+i}(|x|)\right)^2+\sum_{j=1}^{n-t}\left(M_{t+j}(|y|)\right)^2,
\end{eqnarray}
where the third equality follows from (\ref{P-fml}). For any $u=(u^1,u^2)\in F$, we must have $\Vert u^1\Vert_0\le s$ and $\Vert u^2\Vert_0\le t$, which imply that $I_1\setminus{\rm supp}(u^1)\ge m-s$ and $I_2\setminus{\rm supp}(u^2)\ge n-t$. Moreover, it follows from Definition \ref{psiz} that  $\psi_z(u)=(\psi_z^1(u),\psi_z^2(u))$ with $\psi_z^1(u)=x_{\gamma_u^z}u^1$ and $\psi_z^2(u)=\frac{u^2}{x_{\gamma_u^z}}$, where $\gamma_u^z$ is defined by (\ref{gammamin}). Thus,
\begin{eqnarray}\label{case1-nec}
\Vert\psi_z(u)-z\Vert^2&=&\sum_{i\in I_1}|u_ix_{\gamma_u^z}-z_i|^2+\sum_{j\in I_2}|u_j/x_{\gamma_u^z}-z_j|^2\nonumber\\
&=&\sum_{i\in{\rm supp}(u_1)}|u_ix_{\gamma_u^z}-z_i|^2+\sum_{j\in{\rm supp}(u_2)}|u_j/x_{\gamma_u^z}-z_i|^2\nonumber\\
& &+\sum_{i\in I_1\setminus{\rm supp}(u_1)}|z_i|^2+\sum_{j\in I_2\setminus{\rm supp}(u_2)}|z_j|^2\nonumber\\
&\ge&\sum_{i\in I_1\setminus{\rm supp}(u_1)}|z_i|^2+\sum_{j\in I_2\setminus{\rm supp}(u_2)}|z_j|^2\nonumber\\
&\ge&\sum_{i=1}^{m-s}\left(M_{s+i}(|x|)\right)^2+\sum_{j=1}^{n-t}\left(M_{t+j}(|y|)\right)^2\nonumber\\
&=&\Vert \psi_z(u^*)-z\Vert^2,
\end{eqnarray}
where the last equality follows from (\ref{ustar-z}).
Thus, by the definition of like-projection (i.e., Definition \ref{def-proj}), we have $u^*\in\bar{\bf{P}}_F(z)$.

\emph{Case 2}: $z\in\mathbb{R}^{m+n}$ with $x=0$.  In this case, let $u^*=(u^{1*},u^{2*})$ be defined by (\ref{P-fmle1}), then it follows from the definition of $\psi_z$ that $\psi_z(u^*)=(\psi_z^1(u^*),\psi_z^2(u^*))$ with $\psi_z^1(u^*)=u^{1*}$ and $\psi_z^2(u^*)=u^{2*}$. Thus, similar to (\ref{ustar-z}), it follows that
\begin{eqnarray*}
\Vert \psi_z(u^*)-z\Vert^2=1+\sum_{j=1}^{n-t}\left(M_{t+j}(|y|)\right)^2.
\end{eqnarray*}
For any $u=(u^1,u^2)\in F$, we must have $\Vert u^1\Vert_0\le s$ and $\Vert u^2\Vert_0\le t$, which imply that $I_1\setminus{\rm supp}(u^1)\ge n-s$ and $I_2\setminus{\rm supp}(u^2)\ge n-t$. Moreover, it follows from Definition \ref{psiz} that  $\psi_z(u)=(\psi_z^1(u),\psi_z^2(u))$ with $\psi_z^1(u)=u^1$ and $\psi_z^2(u)=u^2$. Thus, in a manner analogous to (\ref{case1-nec}), one obtains
\begin{eqnarray*}
\Vert \psi_z(u)-z\Vert^2&\ge&|1-0|^2+\sum_{j\in I_2\setminus{\rm supp}(u_2)}|z_j|^2\\
&\ge&1+\sum_{j=1}^{n-t}\left(M_{t+j}(|y|)\right)^2\\
&=&\Vert \psi_z(u^*)-z\Vert^2,
\end{eqnarray*}
where the first inequality follows from $u\in F$ and $x=0$. Thus, by the definition of like-projection (i.e., Definition \ref{def-proj}), we have $u^*\in\bar{\bf{P}}_F(z)$.
\ep
\par
The following example reveals the difference between projection and like-projection.
\begin{example}
Suppose that $F$ is defined by (\ref{sprbl-set}), where $m=3,\ n=3,\ s=2,\ t=2$.
\begin{itemize}
\item Given $z=(0,0,3,3,-4,2)^\top$. Then  $x=(0,0,3)^\top$ and $y=(3,-4,2)^\top$. In this case, it is easy to check that
    $$
    {\bf P}_F(z)  = \left\{(1,0,3,3,-4,0)^\top,(0,1,3,3,-4,0)^\top\right\},$$
    $$ {\bar {\bf P}}_F(z) = \{(0,0,1,9,-12,0)^\top\}.
    $$
\item Given $z=(0,0,1,3,-4,2)^\top$. Then $x=(0,0,1)^\top $ and $ y=(3,-4,2)^\top$. In this case, it is easy to check that
    $$
    {\bf P}_F(z)={\bar {\bf P}}_F(z)=\{(0,0,1,3,-4,0)^\top\}.
    $$
\item Given $z=(0,0,0,3,-4,2)^\top$. Then $x=(0,0,0)^\top$ and $ y=(3,-4,2)^\top$. In this case, it is easy to check that
    $$
    {\bar {\bf P}}_F(z)={\bf P}_F(z) =\left\{(1,0,0,0,3,-4,0)^\top,(0,1,0,0,3,-4,0)^\top,(0,0,1,0,3,-4,0)^\top\right\}.
    $$
\end{itemize}
\end{example}
Now, we introduce the definition of an $L$-like stationary point of (\ref{originalmodel}).
\begin{definition}\label{def-L-stp}
Let $F$ be defined by (\ref{sprbl-set}) and let $L>0$ be a given constant.
$z^*=(x^*,y^*)\in F$ is an L-like stationary point of (\ref{originalmodel}) if and only if
$$
z^*\in{\bar {\bf P}}_F(z^*-\nabla f(z^*)/L)\ and\  \nabla_{\gamma_{x^*}}f(z^*)=0,
$$
where $\gamma_{x^*}$ is the first index of the nonzero entries of $z^*$ and ${\bar {\bf P}}_F$ is defined in Definition \ref{def-proj}.
\end{definition}

\begin{remark}
It is worth pointing out the difference between the L-stationary point  and the L-like stationary point.
\begin{itemize}
\item An L-like stationary point of  (\ref{originalmodel}) must be an L-stationary point of (\ref{originalmodel}). In fact, if $z^*$ is an L-like stationary point of  (\ref{originalmodel}), then it follows from Definition \ref{def-L-stp} that $z^*\in F$ and $\nabla_{\gamma_{x^*}}f(z^*)=0$. Since $z^*\in F$, we have $z_{\gamma_{x^*}}^*=1$. This together with $\nabla_{\gamma_{x^*}}f(z^*)=0$ implies that $(z^*-\nabla f(z^*)/L)_{\gamma_{x^*}}=1$. In this case, it follows from Remark \ref{remark3.1} that ${\bar {\bf P}}_F(z^*-\nabla f(z^*)/L)={\bf P}_F(z^*-\nabla f(z^*)/L)$. Thus $z^*$ is an L-stationary point of (\ref{originalmodel}).
\item An L-stationary point of (\ref{originalmodel}) need not be an L-like stationary point of the problem. Such a result holds since the projection of a point $z\in \mathbb{R}^{m+n}$ onto $F$ does not need to be a like-projection of $z$ onto $F$. In addition, if $z^*$ is an L-stationary point of  (\ref{originalmodel}) and $(z^*-\nabla f(z^*)/L)_{\gamma_{x^*}}=1$, then $z^*$ is an L-like stationary point of the problem. In fact, on one hand, it follows from Remark \ref{remark3.1} that ${\bar {\bf P}}_F(z^*-\nabla f(z^*)/L)={\bf P}_F(z^*-\nabla f(z^*)/L)$, and on the other hand, by $(z^*-\nabla f(z^*)/L)_{\gamma_{x^*}}=1$ and $z_{\gamma_{x^*}}^*=1$,
    it follows that $\nabla_{\gamma_{x^*}}f(z^*)=0$. Thus, by Definition \ref{def-L-stp}, $z^*$ is an L-like stationary point of the problem.
\end{itemize}
\end{remark}

Similar to the single-block variable case \cite{BE2013}, if $z^*$ is an $L$-like stationary point of  (\ref{originalmodel}), then the entries of $z^*$ can provide upper bounds for entries of $|\nabla f(z^*)|$. We discuss this fact in the next theorem.

\begin{theorem}
For $z^*=(x^*,y^*)\in F$,  let $\Gamma^*, \Gamma_1^*,\Gamma_2^*$ be given as (\ref{gm-star}) and denote $\gamma_{x^*}$ by the index of the first nonzero entry of $z^*$. Then $z^*$ is an L-like stationary point of (\ref{originalmodel}) if and only if $\nabla f(z^*)\in \mathbb{R}^{m+n}$ is given by
\begin{equation}\label{LSP}
|\nabla_if(z^*)|\left\{
\begin{array}{ll}
=0&\ \mbox{\rm if}\ i\in\Gamma^*,\\
\le LM_s(|x^*|)&\ \mbox{\rm if}\ i\in I_1\setminus\Gamma_1^*,\\
\le LM_t(|y^*|)&\ \mbox{\rm if}\ i\in I_2\setminus\Gamma_2^*,
\end{array}
\right.
\end{equation}
where L is given as in Definition \ref{def-L-stp}.
\end{theorem}

\noindent{\bf Proof}. ($\Rightarrow$) Since $z^*$ is an L-like stationary point of (\ref{originalmodel}), it follows from Definition \ref{def-L-stp} that $z^*\in{\bar {\bf P}}_F(z^*-\nabla f(z^*)/L)$ and $\nabla_{\gamma_{x^*}} f(z^*)=0$. Next, we show that (\ref{LSP}) holds.  Denote
    $$
    \Lambda_1^* := \{i\in I_1: |(z^*-\nabla f(z^*)/L)_i|=M_{s}(|x^*-\nabla_{I_1} f(z^*)/L|)\}, $$
   $$ \Lambda_2^* := \{i\in I_2: |(z^*-\nabla f(z^*)/L)_i|=M_{t}(|y^*-\nabla_{I_2} f(z^*)/L|)\}.$$
    Since $z^*\in{\bar {\bf P}}_F(z^*-\nabla f(z^*)/L)$, we have $x^*\neq0$. Thus, it follows from (\ref{P-fml}) that there exist $\Upsilon_1^*\subseteq\Lambda_1^*$ and $\Upsilon_2^*\subseteq\Lambda_2^*$ such that
    \begin{equation}\label{P-fml-1}
    z_i^*=\left\{
    \begin{array}{ll}
    \frac{(z^*-\nabla f(z^*)/L)_i}{(z^*-\nabla f(z^*)/L)_{\gamma_{\min}^*}}&\ \mbox{\rm if}\; i\in \Omega_1^*,\vspace{1mm}\\
    (z^*-\nabla f(z^*)/L)_{\gamma_{\min}^*}(z^*-\nabla f(z^*)/L)_i&\ \mbox{\rm if}\; i\in \Omega_2^*,\vspace{1mm}\\
    0&\ \mbox{\rm otherwise},
    \end{array}
    \right.
    \end{equation}
   where
    \begin{equation}\label{P-fml-2}
    \left\{\begin{array}{rcl}
    \Omega_1^*&:=&\{i\in I_1: |(z^*-\nabla f(z^*)/L)_i|>M_{s}(|x^*-\nabla_{I_1} f(z^*)/L|)\}\bigcup\Upsilon_1^*,\\
    \\
    \Omega_2^*&:=&\{i\in I_2: |(z^*-\nabla f(z^*)/L)_i|>M_{t}(|y^*-\nabla_{I_2} f(z^*)/L|)\}\bigcup\Upsilon_2^*
    \end{array} \right.
    \end{equation}
    with $\ |\Omega^*_1|=s$ and $|\Omega_2^*|=t$, and
    \begin{eqnarray}\label{P-fml-2-1}
    \gamma_{\min}^*:=\min\ \mbox{\rm supp}(\{(z^*-\nabla f(z^*)/L)_i: i\in \Omega_1^*\}).
    \end{eqnarray}
Note that (\ref{P-fml-1}) is well-defined. In fact, since $z^*\in F$, we have $z^*_{\gamma_{x^*}}=1$, which together with $\nabla_{\gamma_{x^*}}f(z^*)=0$ implies that $(z^*-\nabla f(z^*)/L)_{\gamma_{x^*}}=1$. Thus  ${\rm supp}(\{(z^*-\nabla f(z^*)/L)_i: i\in \Omega_1^*\})$ is non-empty, and hence $(z^*-\nabla f(z^*)/L)_{\gamma_{\min}^*}\neq 0$. So (\ref{P-fml-1}) is well-defined.
Moreover, it is easy to see from (\ref{P-fml-1}) and (\ref{P-fml-2}) that $\mbox{\rm supp}(x^*)=\mbox{\rm supp}(\{(z^*-\nabla f(z^*)/L)_i: i\in \Omega_1^*\})$ and $\mbox{\rm supp}(y^*)=\mbox{\rm supp}(\{(z^*-\nabla f(z^*)/L)_i: i\in \Omega_2^*\})$. Thus, $\gamma_{x^*}=\gamma_{\min}^*$.
    \begin{itemize}
    \item[(a)] We show that $\nabla_if(z^*)=0$ for any $i\in\Gamma^*$. In fact, for any $ i \in \Gamma ^*, $ there are only two cases: $ i \in \Gamma^*_1$ and $  i \in \Gamma^*_2.$
   For the first case, i.e., $i\in\Gamma_1^*$, one has $({\bar {\bf P}}_F(z^*-\nabla f(z^*)/L))_i=z^*_i\neq0$. This, together with (\ref{P-fml-1}), $x^*_{\gamma_{x^*}}=1$ and $\nabla_{\gamma_{x^*}}f(z^*)=0$, implies that
    \begin{eqnarray}\label{thm3-13e1}
    z^*_i-\nabla_if(z^*)/L=z^*_i\left(x^*_{\gamma_{x^*}}-\nabla_{\gamma_x}f(z^*)/L\right) =z^*_i,\;\; \forall i\in\Gamma_1^*.
    \end{eqnarray}
    For the second case, i.e., $i\in\Gamma_2^*$, one has $({\bar {\bf P}}_F(z^*-\nabla f(z^*)/L))_i=z^*_i\neq0$. This, together with (\ref{P-fml-1}), $x^*_{\gamma_{x^*}}=1$ and $\nabla_{\gamma_{x^*}}f(z^*)=0$, implies that
    \begin{eqnarray}\label{thm3-13e2}
    z^*_i-\nabla_if(z^*)/L=z^*_i/\left(x^*_{\gamma_{x^*}}-\nabla_{\gamma_x}f(z^*)/L\right) =z^*_i,\;\; \forall i\in\Gamma_2^*.
    \end{eqnarray}
    Thus, by (\ref{thm3-13e1}) and (\ref{thm3-13e2}), we obtain that $\nabla_if(z^*)=0$ for any $i\in\Gamma^*$.

\item[(b)] We show that $|\nabla_if(z^*)|\leq LM_s(|x^*|)$ for any $i\in I_1\setminus\Gamma_1^*$. There are only the following two cases.

  \emph{Case 1:}  $|\Gamma^*_1|=s$. In this case, it follows from $\mbox{\rm supp}(x^*)=\mbox{\rm supp}(\{(z^*-\nabla f(z^*)/L)_i: i\in \Omega_1^*\})$ that $\mbox{\rm supp}(x^*)\subseteq\Omega^*_1$. Since $|\Gamma^*_1|=|\Omega^*_1|=s$, we have $\Gamma^*_1=\Omega^*_1$. By the definition of $\Omega^*_1$ in (\ref{P-fml-2}), we have that for any $i\in I_1\setminus\Omega^*_1$ and $j\in\Omega^*_1$,
\begin{eqnarray*}
|(z^*-\nabla f(z^*)/L)_i|\le M_s(|x^*-\nabla_{I_1}f(z^*)/L|)\le|(z^*-\nabla f(z^*)/L)_j|.
\end{eqnarray*}
Moreover, it follows from (a) that $\nabla_if(z^*)=0$ for any $i\in\Gamma^*_1$. Thus, combining with $\Gamma^*_1=\Omega^*_1$, it implies that for any $i\in I_1\setminus\Gamma^*_1$ and $j\in\Gamma^*_1$,
\begin{eqnarray*}
|-\nabla_if(z^*)/L|=|(z^*-\nabla f(z^*)/L)_i|\le|(z^*-\nabla f(z^*)/L)_j|=|z^*_j|.
\end{eqnarray*}
Furthermore, it follows that $|-\nabla_if(z^*)/L|\le\min_{j\in\Gamma^*_1} |z^*_j|=M_s(|x^*|)$ for any $i\in I_1\setminus\Gamma^*_1$. Thus, $|\nabla_if(z^*)|\le LM_s(|x^*|)$ for any $i\in I_1\setminus\Gamma^*_1$.

\emph{Case 2:}   $|\Gamma^*_1|<s$. In this case, $M_s(|x^*|)=0$. We must have $|\nabla_if(z^*)|=0\le LM_s(|x^*|)$ for any $i\in I_1\setminus\Gamma_1^*$. Denote
\begin{eqnarray}\label{i0-eq}
i_0:={\rm arg}\max_{j\in I_1\setminus\Gamma^*_1}\left\{|\nabla_jf(z^*)|\right\}.
\end{eqnarray}
Suppose $|\nabla_{i_0}f(z^*)|\neq0$. Since $|\Gamma^*_1|<s=|\Omega^*_1|$, we have $\Gamma^*_1\neq\Omega^*_1$. It follows from (\ref{i0-eq}) that for any $i\in I_1\setminus\Gamma_1^*$,
$$
|(z^*-\nabla f(z^*)/L)_{i_0}|=|\nabla_{i_0} f(z^*)/L|\ge|\nabla_i f(z^*)/L|=|(z^*-\nabla f(z^*)/L)_i|.
$$
Thus $i_0\in\Omega^*_1$ by (\ref{P-fml-2}). For any $u\in\bar{{\bf P}}(z^*-\nabla f(z^*)/L)$, by (\ref{P-fml-1}) and $\gamma_{x^*}=\gamma_{\min}^*$, we have
$$u_{i_0}=\frac{(z^*-\nabla f(z^*)/L)_{i_0}}{z^*_{\gamma_{x^*}}-\nabla_{\gamma_{x^*}}f(z^*)/L}=(z^*-\nabla f(z^*)/L)_{i_0}\neq z^*_{i_0},
$$
 which yields a contradiction. Thus, $|\nabla_if(z^*)|=0$ for any $i\in I_1\setminus\Gamma_1^*$, which implies that $|\nabla_if(z^*)|\leq LM_s(|x^*|)$ for any $i\in I_1\setminus\Gamma_1^*$.

Combining the two cases above yileds $|\nabla_if(z^*)|\leq LM_s(|x^*|)$ for any $i\in I_1\setminus\Gamma_1^*$.

\item[(c)] We show that $|\nabla_if(z^*)|\leq LM_t(|y^*|)$ for any $i\in I_2\setminus\Gamma_2^*$. This can be proved similarly to (b), and thus the proof is omitted.
\end{itemize}
Merging the cases (a)-(c) above  leads to  (\ref{LSP}).

($\Leftarrow$) Suppose that (\ref{LSP}) holds. We show that $z^*$ is an L-like stationary point of   (\ref{originalmodel}). From (\ref{LSP}), one has
\begin{equation}\label{simplified-z-nab}
(z^*-\nabla f(z^*)/L)_i=\left\{
\begin{array}{ll}
z^*_i&\ \mbox{\rm if}\; i\in\Gamma^*,\\
-\nabla_if(z^*)/L&\ \mbox{\rm otherwise.}
\end{array}
\right.
\end{equation}
Since $z^*\in F$, we have $x_{\gamma_{x^*}}^*=1$. This, together with the first equality in (\ref{simplified-z-nab}), implies that $\nabla_{x_{\gamma_{x^*}}^*} f(z^*)=0$. By Definition \ref{def-L-stp}, we need to show that $z^*\in{\bar {\bf P}}_F(z^*-\nabla f(z^*)/L)$.  By Lemma \ref{lem3-4}, we only need to prove there exist $\Omega_1^*$ and $\Omega_2^*$  defined by (\ref{P-fml-2}) such that (\ref{P-fml-1}) holds. To this goal, we divide the remaining proof into  two parts.

{\bf Part 1}. We show that there exists $\Omega^*_1$  defined as (\ref{P-fml-2}) such that \begin{equation}\label{x-pj}
z^*_i=\left\{
\begin{array}{ll}
\frac{(z^*-\nabla f(z^*)/L)_i}{(z^*-\nabla f(z^*)/L)_{\gamma_{\min}^*}}&\mbox{\rm if}\ i\in\Omega^*_1,\\
0&\mbox{\rm if}\ i\in I_1\setminus\Omega^*_1,
\end{array}
\right.
\end{equation}
where $\gamma_{\min}^*$ is defined by (\ref{P-fml-2-1}).
It is sufficient to consider the following two cases.

\emph{Case 1:} $|\Gamma^*_1|=s$. In this case, we define $\Omega^*_1=\Gamma^*_1$. Firstly, we show that $\Omega^*_1$ is just the one defined by (\ref{P-fml-2}). On one hand, it is obvious that $|\Omega^*_1|=s$. On the other hand, since $z^*_i\ge M_s(|x^*|)$ for $i\in \Gamma^*_1$, it follows from (\ref{LSP}) that
    $$
    |z^*_i|\ge |\nabla_jf(z^*)/L|,\ \forall i\in\Gamma^*_1,\ j\in I_1\setminus\Gamma_1^*.
    $$
    Thus, by (\ref{simplified-z-nab}) we have
     \begin{eqnarray}\label{gamma-x-nab}
     |(z^*-\nabla f(z^*)/L)_i|\ge|(z^*-\nabla f(z^*)/L)_j|,\ \forall i\in\Gamma^*_1,\ \ j\in I_1\setminus\Gamma_1^*.
     \end{eqnarray}
     By $|\Gamma^*_1|=s$ and (\ref{gamma-x-nab}), we have $|(x^*-\nabla_{I_1} f(z^*)/L)_i|\ge M_s(|x^*-\nabla_{I_1}f(z^*)/L|)$ for any $i\in\Omega^*_1$ and $|(x^*-\nabla_{I_1} f(z^*)/L)_j|\le M_s(|x^*-\nabla_{I_1}f(z^*)/L|)$ for any $j\in I_1\setminus\Omega^*_1$. Thus $\Omega^*_1$ is just the one defined by (\ref{P-fml-2}) with $\Upsilon_1^*=\Lambda_1^*\cap\Gamma_1^*$.

     Secondly, we show that (\ref{x-pj}) holds.
     By $|\Omega^*_1|=s$ and (\ref{simplified-z-nab}), we have
    $$
    \gamma_{\min}^*=\min\ \mbox{\rm supp}(\{(z^*-\nabla f(z^*)/L)_i: i\in \Omega_1^*\})=\min\ \mbox{\rm supp}(\{z^*_i: i\in \Gamma_1^*\})=\gamma_{x^*}.
    $$
    By $\nabla_{\gamma_{x^*}}f(z^*)=0$, $z^*_{\gamma_{x^*}}=1$ and $\nabla _if(z^*)=0$ for any $i\in\Omega_1^*$, we see that (\ref{x-pj}) holds.

\emph{Case 2:}  $|\Gamma^*_1|=q<s$. In this case, we have $M_s(|x^*|)=0$. This, together with (\ref{LSP}), implies that $\nabla_if(z^*)=0$ for any $i\in I_1$. Thus, $x^*-\nabla_{I_1}f(z^*)/L=x^*$. Choose $\Upsilon_1\subseteq I_1\setminus\Gamma_1^*$ with $|\Upsilon_1|=s-q$ and define $\Omega_1^*:=\Gamma_1^*\cup\Upsilon_1$. Then, $|\Omega_1^*|=s$. Since $x^*-\nabla_{I_1}f(z^*)/L=x^*$, we have $$|(x^*-\nabla f(z^*)/L)_i|=|x^*_i|\ge M_s(|x^*|)=M_s(|x^*-\nabla_{I_1}f(z^*)/L|)$$ for any $i\in\Omega^*_1$  and $$|(x^*-\nabla f(z^*)/L)_j|=|x^*_j|\le M_s(|x^*|)=M_s(|x^*-\nabla_{I_1}f(z^*)/L|)$$ for any $j\in I_1\setminus\Omega^*_1$. Thus $\Omega^*_1$ is just the one defined by (\ref{P-fml-2}). Moreover, we have
    $$
    \gamma_{\min}^*=\min\ \mbox{\rm supp}(\{(z^*-\nabla f(z^*)/L)_i: i\in \Omega_1^*\})=\min\ \mbox{\rm supp}(\{z^*_i: i\in \Omega_1^*\})=\gamma_{x^*}.
    $$
    Thus, by $x^*_{\gamma_{x^*}}=1$ and $\nabla_if(z^*)=0$ for any $i\in I_1$, we see  that (\ref{x-pj}) holds.

{\bf Part 2}. Similar to the proof in {\bf Part 1} for $x$-block, we can also show that there exists $\Omega^*_2\subseteq I_2$ with the same definition as (\ref{P-fml-2}) such that
\begin{equation}\label{y-pj}
z^*_i=\left\{
\begin{array}{ll}
 (z^*-\nabla f(z^*)/L)_{\gamma_{\min}^*}(z^*-\nabla f(z^*)/L)_i&\mbox{\rm if}\ i\in\Omega^*_2,\\
0&\mbox{\rm if}\ i\in I_2\setminus\Omega^*_2.
\end{array}
\right.
\end{equation}

Combining (\ref{x-pj}) and (\ref{y-pj}) yields (\ref{P-fml-1}), and thus, by Lemma \ref{lem3-4}, we obtain that $z^*\in{\bar {\bf P}}_F(z^*-\nabla f(z^*)/L)$. Furthermore, by Definition \ref{def-L-stp},
$z^*$ is an L-like stationary point of the problem (\ref{originalmodel}).
\ep

The relation between $L$-like stationarity and N-type stationarity is stated  in the theorem below.

\begin{theorem}\label{L-to-B}
If $z^*=(x^*,y^*)$ is an $L$-like stationary point of (\ref{originalmodel}), then $z^*$ is also an ${\rm N^B}$-stationary point of (\ref{originalmodel}).
\end{theorem}

\noindent{\bf Proof}. Consider the following four cases:
\begin{itemize}
\item $\Vert x^*\Vert_0=s$ and $\Vert y^*\Vert_0=t$.
For this case, $M_s(|x^*|)>0$ and $M_t(|y^*|)>0$. Combining with (\ref{LSP}), we see that  $\nabla_if(z^*)=0$ for any $i\in\Gamma^*$, which implies  (\ref{N-B-c1}).

\item  $\Vert x^*\Vert_0<s$ and $\Vert y^*\Vert_0=t$.
In this case, $M_s(|x^*|)=0$ and $M_t(|y^*|)>0$. Combining with (\ref{LSP}) leads to  $\nabla_if(z^*)=0$ for any $i\in I_1\cup\Gamma_2^*$, which implies   (\ref{N-B-c1}).

\item  $\Vert x^*\Vert_0=s$ and $\Vert y^*\Vert_0<t$.
In this case, $M_s(|x^*|)>0$ and $M_t(|y^*|)=0$. Combining with (\ref{LSP}), we have $\nabla_if(z^*)=0$ for any $i\in\Gamma^*_1\cup I_2$, which implies (\ref{N-B-c1}).

\item  $\Vert x\Vert_0<s$ and $\Vert y\Vert_0<t$.
In this case, $M_s(|x^*|)=0$ and $M_t(|y^*|)=0$. Combining with (\ref{LSP}) yields  $\nabla f(z^*)=0$ for any $i\in I$, which also implies  (\ref{N-B-c1}) .
\end{itemize}
  Thus, by Theorem \ref{thm-N-B}, the desired result follows.
\ep

The following two theorems provide necessary optimal conditions for (\ref{originalmodel}) through the L-like stationary points.

\begin{theorem}\label{L-opt1}
Let $z^*=(x^*,y^*)$ be an optimal solution of (\ref{originalmodel}) and  $A:=\mathscr{A}\times_3 y^*$. If $\nabla_{\gamma_x}f(z^*)=0$ and $$L>\max\left\{\lambda_{max}\left(A^\top A\right),\max_{i\in I_2\setminus\Gamma^*_2}\frac{|\nabla_if(z^*)|}{M_t(|y^*|)}\right\},$$
then $z^*$ is an $L$- like stationary point of (\ref{originalmodel}), and $\bar{{\bf P}}_F(z^*-\nabla f(z^*)/L)$ is a singleton.
Thus
\begin{eqnarray}\label{stnry-eqt}
z^*=\bar{{\bf P}}_F\left(z^*-\nabla f(z^*)/L\right).
\end{eqnarray}
\end{theorem}

\noindent{\bf Proof}. Since $z^*=(x^*,y^*)$ is an optimal solution of (\ref{originalmodel}), it follows from Theorem \ref{B-nec} that $z^*$ is an ${\rm N^B}$-stationary point of (\ref{originalmodel}). Thus, by Theorem \ref{thm-N-B}, we have that
$\nabla_if(z^*)=0$ for any $i\in I_2$ if $\Vert y^*\Vert_0<t$, and that
 $\nabla_if(z^*)=0$ for any $i\in\Gamma_2^*$ if $\Vert y^*\Vert_0=t$, where $\Gamma_2^*$ is defined as (\ref{gm-star}). Let $z=(x,y)\in\bar{\bf P}_F(z^*-\nabla f(z^*)/L)$. We show that $z=z^*$.

{\bf Part 1}. We show $y=y^*$. Consider the following two cases.

\emph{ Case 1}:  $\Vert y^*\Vert_0=q<t$. In this case, since $\nabla_if(z^*)=0$ for any $i\in I_2$, we have $y^*-\nabla_{I_2}f(z^*)/L=y^*$. Pick $\Upsilon_2\subseteq I_2\setminus \Gamma_2^*$ with $|\Upsilon_2|=t-q$ and denote $\Omega_2=\Upsilon_2\cup\Gamma^*_2$. Since $\nabla_{I_2}f(z^*)=0$, we have $|(z^*-\nabla f(z^*)/L)_i|=|y^*_i|$ for any $i\in\Omega_2$ and $M_t(|y^*-\nabla_{I_2}f(z^*)/L|)= M_t(|y^*|)$. Then for any $i\in\Omega_2$, it follows from $|y^*_i|\ge M_t(|y^*|)$ that  $|(z^*-\nabla f(z^*)/L)_i|\ge M_t(|y^*-\nabla_{I_2}f(z^*)/L|)$. Thus, $\Omega_2$ satisfies (\ref{P-fml-2}) with $\Upsilon_2$ being chosen arbitrarily. Note that
  $\Upsilon_2$ is picked arbitrarily, which leads to the arbitrarity of $\Omega_2$. It follows from $z=(x,y)\in\bar{{\bf P}}_F\left(z^*-\nabla f(z^*)/L\right)$ and (\ref{P-fml}) that
    \begin{equation}\label{gm1s}
    y_j=\left\{\begin{array}{ll}
    \left(x^*_{\gamma_{x^*}}-\nabla_{\gamma_{x^*}}f(z^*)/L\right)y^*_j&\mbox{\rm if}\ j\in\Omega_2,\\
    0&\mbox{\rm if}\ j\in I_2\setminus\Omega_2.
    \end{array}\right.
    \end{equation}
Since $x^*_{\gamma_{x^*}}=1,\nabla_{\gamma_{x^*}}f(z^*)=0$, we have $y=y^*$.

 \emph{ Case 2}: $\Vert y^*\Vert_0=t$. In this case, we have $\nabla_if(z^*)=0$ for any $i\in\Gamma_2^*$. Then
    \begin{equation*}
    (y^*-\nabla_{I_2}f(z^*)/L)_i=\left\{\begin{array}{ll}
    y^*_i &\mbox{\rm if}\ i\in\Gamma^*_2,\\
   -\nabla_if(z^*)/L &\mbox{\rm if}\ i\in I_2\setminus\Gamma^*_2.
    \end{array}\right.
    \end{equation*}
    As $L>\max_{i\in I_2\setminus\Gamma^*_2}\frac{|\nabla_if(z^*)|}{M_t(|y^*|)}$, we have $|\nabla_if(z^*)|/L<M_t(|y^*|)$ for any $i\in I_2\setminus\Gamma_2^*$. Since $|\Gamma_2^*|=t$, we have $|z^*_i|\ge M_t(y^*)$ for any $i\in\Gamma_2^*$. Thus
    \begin{eqnarray*}
    &&|(y^*-\nabla_{I_2}f(z^*)/L)_{i_1}|\ge\cdots\ge|(y^*-\nabla_{I_2}f(z^*)/L)_{i_t}|\\
    &&>|(y^*-\nabla_{I_2}f(z^*)/L)_{j_1}|\ge\cdots\ge|(y^*-\nabla_{I_2}f(z^*)/L)_{j_{n-t}}|,
    \end{eqnarray*}
    where $i_1,\ldots,i_t\in\Gamma_2^*, $ $ j_1,\ldots,j_{n-t}\in I_2\setminus\Gamma_2^*$ and $i_1,\ldots,i_t,j_1,\ldots,j_{n-t}$ are not equal to each other. Then the only way to choose $\Omega_2$ satisfying (\ref{P-fml-2}) is $\Omega_2=\Gamma^*_2$. Since $\nabla_if(z^*)=0$ for any $i\in\Omega_2$, then we also obtain (\ref{gm1s}). Thus, we must have $y=y^*$.

{\bf Part 2}. We show $x=x^*$. Suppose that there exists $z=(x,y^*)$ with $\Vert x-x^*\Vert>0$ such that $z\in\bar{{\bf P}}_F(z^*-\nabla f(z^*)/L)$.   Combining $z^*_{\gamma_x}=1$ with $\nabla_{\gamma_x}f(z)=0$, One has
$$\Vert z-\left(z^*-\nabla f(z^*)/L\right)\Vert^2\le\Vert z^*-\left(z^*-\nabla f(z^*)/L\right)\Vert^2,$$
which implies that
$$ \left<z-\left(z^*-\nabla f(z^*)/L\right),z-\left(z^*-\nabla f(z^*)/L\right)\right>\le \frac{1}{L^2}\Vert\nabla f(z^*)\Vert^2,$$
 and hence $  \left<z-z^*,\nabla f(z^*)\right>\le-\frac{L}{2}\Vert z-z^*\Vert^2$  which is equivalent to
\begin{equation} \left<x-x^*, \nabla_{I_1} f(z^*)\right>\le-\frac{L}{2}\Vert x-x^*\Vert^2.\label{projection01}
\end{equation}
Let $\phi_{y^*}(x)=\frac{1}{2}\Vert Ax-b\Vert^2$. By expanding $\phi_{y^*}$ at $x^*$ and using (\ref{projection01}) , one has
\begin{align*}
f(z)&=\phi_{y^*}(x)\le\phi_{y^*}(x^*)+\left<x-x^*,\nabla\phi_{y^*}(x^*)\right>+\frac{\lambda_{max}(A^\top A)}{2}\Vert x-x^*\Vert^2\\
&   \le  \phi_{y^*}(x^*)-\frac{L}{2}\Vert x-x^*\Vert^2+\frac{\lambda_{max}(A^\top A)}{2}\Vert x-x^*\Vert^2\\
&<\phi_{y^*}(x^*)=f(z^*),
\end{align*}
which contradicts the optimality of $z^*$. Combining {\bf Part 1} with {\bf Part 2} yields the equality (\ref{stnry-eqt}).
\ep

By a similar proof to the above, we can also show the next result. For completeness, we still give a full proof for this result.
\begin{theorem}
Let $z^*=(x^*,y^*)$ be an optimal solution of (\ref{originalmodel}) and $B:=\mathscr{A}\times_2 x^*$. If $\nabla_{\gamma_x}f(z^*)=0$ and $$L>\max\left\{\lambda_{max}\left(B^\top B\right),\max_{i\in I_1\setminus\Gamma_1^*}\frac{|\nabla_if(z^*)|}{M_s(|x^*|)}\right\},$$
then $z^*$ is an $L$-like stationary point of (\ref{originalmodel}), and $\bar{{\bf P}}_F(z^*-\nabla f(z^*)/L)$ is a singleton.
Thus  (\ref{stnry-eqt}) holds, i.e.,
\begin{eqnarray*}
z^*=\bar{{\bf P}}_F\left(z^*-\nabla f(z^*)/L\right).
\end{eqnarray*}
\end{theorem}

\noindent{\bf Proof}. While the proof is similar to that of  Theorem \ref{L-opt1}. We still provide the proof here for completeness. Since $z^*=(x^*,y^*)$ is an optimal solution of (\ref{originalmodel}), it follows from Theorem \ref{B-nec} that $z^*$ is an ${\rm N^B}$-stationary point of (\ref{originalmodel}). Thus, by Theorem \ref{thm-N-B}, we have that  $\nabla_if(z^*)=0$ for any $i\in I_1$ if $\Vert x^*\Vert_0<s$, and that
 $\nabla_if(z^*)=0$ for any $i\in\Gamma_1^*$ if $\Vert x^*\Vert_0=s$, where $\Gamma_1^*$ is defined as (\ref{gm-star}). Let $z=(x,y)\in\bar{\bf P}_F(z^*-\nabla f(z^*)/L)$. We show that $z=z^*$.

{\bf Part 1}. We show $x=x^*$. Consider the   two cases below.

\emph{Case 1:}   $\Vert x^*\Vert_0=q<s$. In this case, since $\nabla_if(z^*)=0$ for any $i\in I_1$, we have $x^*-\nabla_{I_1}f(z^*)/L=x^*$. Pick $\Upsilon_1\subseteq I_1\setminus \Gamma_1^*$ with $|\Upsilon_1|=t-q$ and denote $\Omega_1=\Upsilon_1\cup\Gamma^*_1$. Since $\nabla_{I_1}f(z^*)=0$, we have $|(z^*-\nabla f(z^*)/L)_i|=|x^*_i|$ for any $i\in\Omega_1$ and $M_s(|x^*-\nabla_{I_1}f(z^*)/L|)= M_s(|x^*|)$. Then for any $i\in\Omega_1$, it follows from $|x^*_i|\ge M_s(|x^*|)$ that  $|(z^*-\nabla f(z^*)/L)_i|\ge M_s(|x^*-\nabla_{I_1}f(z^*)/L|)$. Thus, $\Omega_1$ satisfies (\ref{P-fml-2}) with $\Upsilon_1$ being chosen arbitrarily. Note that
   $\Upsilon_1$ is picked arbitrarily, which leads to the arbitrarity of $\Omega_1$. It follows from $z=(x,y)\in\bar{{\bf P}}_F\left(z^*-\nabla f(z^*)/L\right)$ and (\ref{P-fml}) that
    \begin{equation}\label{gm2s}
    x_j=\left\{\begin{array}{ll}
    x^*_j/\left(x^*_{\gamma_{x^*}}-\nabla_{\gamma_{x^*}}f(z^*)/L\right)&\mbox{\rm if}\ j\in\Omega_1,\\
    0&\mbox{\rm if}\ j\in I_1\setminus\Omega_1.
    \end{array}\right.
    \end{equation}
Since $x^*_{\gamma_{x^*}}=1,\nabla_{\gamma_{x^*}}f(z^*)=0$, we have $x=x^*$.

\emph{Case 2:}  $\Vert x^*\Vert_0=s$. In this case, we have $\nabla_if(z^*)=0$ for any $i\in\Gamma_1^*$. Then we have
    \begin{equation*}
    (x^*-\nabla_{I_1}f(z^*)/L)_i=\left\{\begin{array}{ll}
    x^*_i &\mbox{\rm if}\ i\in\Gamma^*_1,\\
    -\nabla_if(z^*)/L &\mbox{\rm if}\ i\in I_1\setminus\Gamma^*_1.
    \end{array}\right.
    \end{equation*}
    Since $L>\max_{i\in I_1\setminus\Gamma^*_1}\frac{|\nabla_if(z^*)|}{M_s(|x^*|)}$, we have $|\nabla_if(z^*)|/L<M_s(|x^*|)$ for any $i\in I_1\setminus\Gamma_1^*$. Since $|\Gamma_1^*|=s$, we have $|z^*_i|\ge M_s(x^*)$ for any $i\in\Gamma_1^*$. We have
    \begin{eqnarray*}
    &&|(x^*-\nabla_{I_1}f(z^*)/L)_{i_1}|\ge\cdots\ge|(x^*-\nabla_{I_1}f(z^*)/L)_{i_t}|\\
    &&>|(x^*-\nabla_{I_1}f(z^*)/L)_{j_1}|\ge\cdots\ge|(x^*-\nabla_{I_1}f(z^*)/L)_{j_{n-t}}|,
    \end{eqnarray*}
    where $i_1,\ldots,i_t\in\Gamma_1^*,$ $ j_1,\ldots,j_{n-t}\in I_1\setminus\Gamma_1^*$ and $i_1,\ldots,i_t,j_1,\ldots,j_{n-t}$ are not equal to each other. Then the only way to choose $\Omega_1$ satisfying (\ref{P-fml-2}) is $\Omega_1=\Gamma^*_1$. Since $\nabla_if(z^*)=0$ for any $i\in\Omega_1$, then we also have (\ref{gm2s}) holds. Thus, we must have $x=x^*$.

{\bf Part 2}. We show $y=y^*$. Suppose that there exists $z=(x^*,y)$ with $\Vert y-y^*\Vert>0$ such that $z\in\bar{{\bf P}}_F(z^*-\nabla f(z^*)/L)$.   Combine $z^*_{\gamma_x}=1$ with $\nabla_{\gamma_x}f(z)=0$, we have
$$ \Vert z-\left(z^*-\nabla f(z^*)/L\right)\Vert^2\le\Vert z^*-\left(z^*-\nabla f(z^*)/L\right)\Vert^2,$$
which implies that  $$ \left<z-\left(z^*-\nabla f(z^*)/L\right),z-\left(z^*-\nabla f(z^*)/L\right)\right>\le \frac{1}{L^2}\Vert\nabla f(z^*)\Vert^2,$$
and thus $\left<z-z^*,\nabla f(z^*)\right>\le-\frac{L}{2}\Vert z-z^*\Vert^2 $ which is equivalent to
\begin{equation}  \left<y-y^*, \nabla_{I_2} f(z^*)\right>\le-\frac{L}{2}\Vert y-y^*\Vert^2. \label{projection02}
\end{equation}
Let $\phi_{x^*}(y)=\frac{1}{2}\Vert By-b\Vert^2$. By expanding $\phi_{x^*}$ at $y^*$ and (\ref{projection02}), one has
\begin{align*}
f(z)&=\phi_{x^*}(y)\le\phi_{x^*}(y^*)+\left<y-y^*,\nabla\phi_{x^*}(y^*)\right>+\frac{\lambda_{max}(B^\top B)}{2}\Vert y-y^*\Vert^2\\
& \le \phi_{x^*}(y^*)-\frac{L}{2}\Vert y-y^*\Vert^2+\frac{\lambda_{max}(B^\top B)}{2}\Vert y-y^*\Vert^2\\
&< \phi_{x^*}(y^*)=f(z^*),
\end{align*}
which contradicts  the optimality of $z^*$. Combining {\bf Part 1} with {\bf Part 2} leads to  (\ref{stnry-eqt}).
\ep

\vskip 0.08in

For problem (\ref{originalmodel}) and fixed $z^*=(x^*,y^*)\in\mathbb{R}^{m+n}$, denote $A:=\mathscr{A}\times_2x^*$ and $B:=\mathscr{A}\times_3y^*$. Let $\phi_{y^*}(x)=\frac{1}{2}\Vert Bx-b\Vert^2$ and $\phi_{x^*}(y)=\frac{1}{2}\Vert Ay-b\Vert^2$ . Then we obtain the two subproblems below:
\begin{equation}\label{subproblem-2}
\min \{ \phi_{y^*}(x)=\frac{1}{2}\Vert Bx-b\Vert^2 :
  ~ \|x\|_0\le s\}.
\end{equation}
\begin{equation}\label{subproblem-1}
\min \{ \phi_{x^*}(y)=\frac{1}{2}\Vert Ay-b\Vert^2:
~\|y\|_0\le t\},
\end{equation}

\begin{theorem}
For some $L>0$, if $z^*:=(x^*,y^*)$ is an $L$-like stationary point of (\ref{originalmodel}), then $x^*$ and $y^*$ are the $L$-stationary points of (\ref{subproblem-2}) and (\ref{subproblem-1}), respectively. Conversely, if $x^*$ and $y^*$ are $L$-stationary points of (\ref{subproblem-2}) and (\ref{subproblem-1}), respectively, then $\bar{z}^*:=(\bar{z}^*_1,\bar{z}^*_2)$ is an $N^C$-stationary point of (\ref{originalmodel}), where $\bar{z}^*_1:=x^*/x^*_{\gamma_x}$ and $\bar{z}^*_1:=x^*_{\gamma_x}y^*$.
\end{theorem}

\noindent{\bf Proof}. Since $z^*:=(x^*,y^*)$ is an $L$-like stationary point of (\ref{originalmodel}), combining (\ref{LSP}) with $\nabla f(z^*)=(\nabla\phi_{y^*}(x^*),\nabla\phi_{x^*}(y^*))$, we obtain that
\begin{equation}\label{x-L-stp}
|\nabla_i\phi_{y^*}(x^*)|\left\{
\begin{array}{ll}
=0&\mbox{\rm if}\ i\in\Gamma_1^*,\\
\le LM_s(|x^*|)&\mbox{\rm if}\ i\in I_1\setminus\Gamma_1^*
\end{array}
\right.
\end{equation}
and
\begin{equation}\label{y-L-stp}
|\nabla_i\phi_{x^*}(y^*)|\left\{
\begin{array}{ll}
=0&\mbox{\rm if}\ i\in\Gamma_2^*,\\
\le LM_t(|y^*|)&\mbox{\rm if}\ i\in I_2\setminus\Gamma_2^*,
\end{array}
\right.
\end{equation}
respectively. It follows from Lemma 2.2 in \cite{BE2013} that $x^*$ and $y^*$ are $L$-stationary points of (\ref{subproblem-2}) and (\ref{subproblem-1}), respectively.

Conversely, since $x^*$ and $y^*$ are $L$-stationary points of (\ref{subproblem-2}) and (\ref{subproblem-1}), respectively, we see that (\ref{x-L-stp}) and (\ref{y-L-stp}) hold. Thus, $\nabla_i\phi_{y^*}(x^*)=0$ for any $i\in\Gamma^*_1$ and $\nabla_j\phi_{x^*}(y^*)=0$ for any $j\in\Gamma^*_2$. Furthermore, we have $\nabla_if(\bar{z}^*)=\nabla_i\phi_{\bar{z}^*_2}(\bar{z}^*_1)$  for any $i\in\Gamma^*_1$, where
$$
\phi_{\bar{z}^*_2}(\bar{z}^*_1)=\frac{1}{2}\Vert(\mathscr{A}\times_3\bar{z}^*_2)\bar{z}^*_1-b\Vert^2=\frac{1}{2}\Vert\left(\mathscr{A}\times_3(x^*_{\gamma_{x^*}}y^*)\right)(x^*/x^*_{\gamma_x})-b\Vert^2.
$$
So, for any $i\in\Gamma^*_1$, we have
$$
\nabla_if(\bar{z}^*) = \left(\mathscr{A}\times_3(x^*_{\gamma_x}y^*)^\top(\mathscr{A}x^*y^*-b)\right)_i
 = x^*_{\gamma_x}\left(B^\top(Bx^*-b)\right)_i
 = x^*_{\gamma_x}\nabla_i\phi_{y^*}(x^*)=0.
$$
Symmetrically, we have
\begin{eqnarray*}
\nabla_jf(z)=\nabla_i\phi_{x^*}(y^*)/x^*_{\gamma_x}=0,\ \forall j\in\Gamma^*_2.
\end{eqnarray*}
By Theorem \ref{thm-N-C}, we conclude that $z^*$ is an $N^C$-stationary point of (\ref{originalmodel}).
\ep

Recall that for the sparsity constrained optimization (\ref{nonlinear-model}), as shown in \cite{PXZ2015}, the concepts of $L$-stationarity and $N^B$-stationarity are equivalent when the sparsity constraint is not active. However, similar result for $L$-like stationarity and $N^B$-stationarity does not hold. In fact, $L$-like stationarity remains stronger than $N^B$-stationarity for the problem (\ref{originalmodel}) even when the sparsity constraint is not active. We present the following example to illustrate it.

\begin{example}
Consider the problem (\ref{originalmodel}) with $b=(1,7)^\top$, tensor $\mathscr{A}=(a_{j_1j_2j_3})\in\mathbb{R}^{2\times4\times4}$, where $a_{111}=1,a_{122}=1,a_{123}=1,a_{131}=1,a_{134}=1,a_{211}=1,a_{213}=2,a_{221}=1,a_{222}=2,a_{232}=4,a_{244}=1$ and all other entries are zero. Then for any $z=(x,y)\in\mathbb{R}^8$, the objective function in (\ref{originalmodel}) is defined by $f(z):=\frac{1}{2}\Vert \mathscr{A}xy-b\Vert^2$ with
\begin{equation*}
\mathscr{A}xy-b=\left[
\begin{array}{l}
x_1y_1+x_2y_2+x_2y_3+x_3y_1+x_3y_4-1\\
x_1y_1+2x_1y_3+x_2y_1+2x_2y_2+4x_3y_2+x_4y_4-7
\end{array}
\right].
\end{equation*}
Let $s=t=3$. For $z^*=(1,1,0,0,2,1,0,0)^\top$, we have $\gamma_{x^*}=1,z^*_{\gamma_{x^*}}=1,|\Gamma^*_1|=2<s,|\Gamma^*_2|=2<t$, then $z^*\in F$. It is straightforward to verify that $\nabla f(z^*)=(2,0,0,0,0,0,0,0)^\top$, i.e., $\nabla_if(z^*)=0$ for any $i>\gamma_{x^*}$, which implies that $z^*$ is an $N^B$-stationary point. However, for any $L\in(0,+\infty)$,
$$
z^*-\nabla f(z^*)/L=(1-2/L,1,0,0,2,2,0,0)^\top
$$
and
$$
{\bf \bar{P}}_F(z^*-\nabla f(z^*)/L)=\{(1,1/(1-2/L),0,0,2-4/L,2-4/L,0,0)\},
$$
which imply $z^*\notin{\bf \bar{P}}_F(z^*-\nabla f(z^*)/L)$. Thus, $z^*$ is not an $L$-like stationary point.
\end{example}

\subsection{M-stationary point}

Recall that Burdakov, Kanzow and Schwartz \cite{BKS2016} reformulated and relaxed the single-block-variable optimization problem with sparse constraint as a smooth model with complementarity constraint. We also consider complementarity-type reformulation of  (\ref{originalmodel}), based on which we introduce the M-stationary point of (\ref{originalmodel}). For this purpose, we introduce the auxiliary variable $w=(u,v)\in\mathbb{R}^{m+n}$ and consider the following mixed-integer model:
\begin{eqnarray}\label{mix-int-model}
\begin{array}{cl}
\min & f(x,y)=\frac{1}{2}\|\mathscr{A} xy-b\|^2\vspace{1mm}\\
\mbox{\rm s.t.} & e_{\gamma_x}^\top x=1,\vspace{1mm}\\
& z^\top w=0,\ e_{I_1}^\top u\ge m-s,\ e_{I_2}^\top v\ge n-t,\vspace{1mm}\\
&w_i\in\left\{0,1\right\},\forall i\in I.
\end{array}
\end{eqnarray}
The following theorem shows the relationship between the global minima of the problems (\ref{originalmodel}) and (\ref{mix-int-model}).

\begin{theorem}\label{optm-int-mix}
$z^*\in\mathbb{R}^{m+n}$ is an optimal solution of  (\ref{originalmodel}) if and only if there exists $w^*\in\mathbb{R}^{m+n}$ such that the pair $(z^*,w^*)$ is an optimal solution of  (\ref{mix-int-model}).
\end{theorem}

\noindent{\bf Proof}. ($\Rightarrow$) Suppose that $z^*\in\mathbb{R}^{m+n}$ is an optimal solution of  (\ref{originalmodel}). Then, we have $z^*\in F$. Let $w^*=(w^{1*},w^{2*})\in \mathbb{R}^{m+n}$ be given as
\begin{equation*}
w^*_i=\left\{
\begin{array}{ll}
0&\mbox{\rm if}\ z^*_i\neq0,\\
1&\mbox{\rm if}\ z^*_i=0,
\end{array}\;\; \forall i\in I.
\right.
\end{equation*}
Then ${z^*}^\top w^*=0$, $e_{I_1}^\top w^{1*}\ge n-s$ and $e_{I_2}^\top w^{2*}\ge n-t$. Thus $(z^*,w^*)$ is a feasible point of  (\ref{mix-int-model}). Suppose that $(z^*,w^*)$ is not an optimal solution of  (\ref{mix-int-model}), then there exists a feasible point $(\bar{z},\bar{w})$ of  (\ref{mix-int-model}) with $\bar{z}=(\bar{x},\bar{y})$ and $\bar{w}=(\bar{u},\bar{v})$ such that $f(\bar{z})<f(z^*)$. Since $\bar{z}^\top\bar{w}=0$, $e_{I_1}^\top\bar{u}\ge m-s$, $e_{I_2}^\top\bar{v}\ge n-t$ and $\bar{w}_i\in\left\{0,1\right\}$ for all $i\in I$, it follows that $\|\bar{x}\|_0\leq s$ and $\|\bar{y}\|_0\leq t$. These together with $e_{\gamma_{\bar{x}}}^\top \bar{x}=1$ imply that $\bar{z}\in F$. However, $f(\bar{z})<f(z^*)$, which contradicts the assumption that $z^*$ is an optimal solution of  (\ref{originalmodel}).

($\Leftarrow$)  Suppose that there exists $w^*\in\mathbb{R}^{m+n}$ such that $(z^*,w^*)$ is an optimal solution of (\ref{mix-int-model}). Then one must have $z^*\in F$. Suppose that $z^*$ is not an optimal solution of  (\ref{originalmodel}). Then there exists $\bar{z}\in F$ such that $f(\bar{z})<f(z^*)$. Choose $\bar{w}\in\mathbb{R}^{m+n}$ such that
\begin{equation*}
\bar{w}_i=\left\{
\begin{array}{ll}
0&\mbox{\rm if}\ \bar{z}_i\neq0,\\
1&\mbox{\rm if}\ \bar{z}_i=0,
\end{array}\;\; \forall i\in I.
\right.
\end{equation*}
Then $(\bar{z},\bar{w})$ is a feasible point of problem (\ref{mix-int-model}). However, $f(\bar{z})<f(z^*)$. These contradict the assumption that $(z^*,w^*)$ is an optimal solution to problem (\ref{mix-int-model}).
\ep

In order to reduce the complexity of $0$-$1$ constraints in problem (\ref{mix-int-model}), we consider the following relaxed model.
\begin{eqnarray}\label{relaxed-model}
\begin{array}{cl}
\min & f(x,y)=\frac{1}{2}\|\mathscr{A} xy-b\|^2\vspace{1mm}\\
\mbox{\rm s.t.} & e_{\gamma_x}^\top x=1,\vspace{1mm}\\
& z^\top w=0, e_{I_1}^\top u\ge m-s, e_{I_2}^\top v\ge n-t,\vspace{1mm}\\
&w_i\in[0,1], \forall i\in I.
\end{array}
\end{eqnarray}

Problems (\ref{originalmodel}) and  (\ref{relaxed-model}) share common global minima, as indicated by the following theorem, to which the proof is similar to that of Theorem  \ref{optm-int-mix} and thus omitted.

\begin{theorem}
$z^*\in\mathbb{R}^{m+n}$ is an optimal solution of    (\ref{originalmodel}) if and only if there exists $w^*\in\mathbb{R}^{m+n}$ such that $(z^*,w^*)$ is an optimal solution of  (\ref{relaxed-model}).
\end{theorem}

Denote
$
E:=\left\{z=(x,y)
\in(\mathbb{R}^{m}\setminus\left\{0\right\})\times\mathbb{R}^{n}:\ x_{\gamma_x}=1\right\}, $ and
$$ g(w):=\left(
\begin{array}{c}
(m-s)-e_{I_1}^\top u\\
(n-t)-e_{I_2}^\top v\\
-w\\
w-e_I
\end{array}
\right).
$$
Then, (\ref{relaxed-model}) can be rewritten as
\begin{eqnarray}\label{relaxed-model-1}
\begin{array}{cl}
\min & f(x,y)=\frac{1}{2}\|\mathscr{A} xy-b\|^2\\
\mbox{\rm s.t.} & g(w)\le0,\ z^\top w=0,\\
& z\in E.
\end{array}
\end{eqnarray}
Denote the feasible set of (\ref{relaxed-model-1}) by $\bar{F}$, and denote  the normal cone of $E$ at $z\in E$ by $N^\#_{E}(z)$, where $\#\in\left\{B,C\right\}$ means the Bouligand sense and Clarke sense, respectively. We provide the following lemma, which gives the formula for $N^\#_{E}(z)$.

\begin{lemma}\label{N-E-z}
For any given $z\in E$, the normal cone of $E$ at $z$ is
\begin{eqnarray}\label{NEz}
N^\#_{E}(z)=\left\{d\in\mathbb{R}^{m+n}:\ d_i=0,\ \text{for\ any} \ i>\gamma_x\right\}.
\end{eqnarray}
\end{lemma}

\noindent{\bf Proof}. Denote
$Q:=\left\{d\in\mathbb{R}^{m+n}:\ d_i=0,\ \text{for\ any}\ i\le\gamma_x\right\}$. Then, by Definitions \ref{Bcone} and \ref{Ccone}, we need to show
$T^\#_E(z)=Q$, where $T^\#_E(z)$ denotes the tangent cone of $E$ at $z$.

 (i) We  prove  that $T^B_E(z)=Q$. Firstly, we show $T^B_E(z)\subseteq Q$. For any $d\in T^B_E(z)$, it follows from Definition \ref{Bcone} that there exist $\left\{z^k\right\}\subset E$ satisfying ${\rm lim}_{k\to\infty}z^k=z$ and $\lambda_k\in\mathbb{R}^+$ such that ${\rm \lim}_{k\to\infty}\lambda_k(z^k-z)=d$. By Lemma \ref{x-k-1}, there exists $k_0\in\mathbb{N}$ such that $z^k_i=0$ for any $i<\gamma_x$ and $z^k_{\gamma_x}=1$ for all $k\geq k_0$. Thus, one has
$$
d_i=\lim_{k\to\infty}\lambda_k(z^k_i-z_i)=0,\ \forall i\le\gamma_x.
$$
That is $d\in Q$. Thus we have $T_E(z)\subseteq Q$. Secondly, we show $T_E(z)\supseteq Q$. For any $d\in Q$, let $z^k=z+\frac{d}{k}$ for all $k\in\mathbb{N}$. By $z\in E$, one has $z^k_i=0$ for any $i<\gamma_x$ and $z^k_{\gamma_x}=1$ for all $k\in\mathbb{N}$. Thus $\left\{z^k\right\}\subset E$. Besides, we also have ${\rm lim}_{k\to\infty}z^k=z$ and ${\rm \lim}_{k\to\infty}\frac{1}{k}(z^k-z)=d$. That is $d\in T^B_E(z)$.  Thus we have $T^B_E(z)\supseteq Q$. Thus $T^B_E(z)=Q$.

(ii) We prove  that $T^C_E(z)=Q$. Firstly, we show $T^C_E(z)\subseteq Q$. Suppose there exits $d\in T^C_E(z)$ with $d_i\neq0$ where $i\le\gamma_x$. Then for any $\left\{d^k\right\}\subseteq E$ with $\lim_{k\to\infty}d^k=d$, there exists $k_1>0$ such that $d^k_i\neq0$ for all $k>k_1$. Without loss of generality, suppose $i$ is the first index of nonzero entries of $d^k$. Pick $\lambda_k:=\frac{1}{k|d^k_i|}$ for all $k\in\mathbb{N}$. For any $z^k\to z$, there exists $k_2>0$ such that $z^k_j=0$ for all $j<\gamma_x$ and $z^k_{\gamma_x}=1$ for all $k>k_2$. It follows that $i$ is the first index of nonzero entries of $z^k+\lambda_kd^k$. Let $K=\max\left\{k_1,k_2\right\}$, one has
\begin{equation*}
z_i^k+\lambda_kd_i^k=\left\{
\begin{array}{ll}
\mbox{\rm sgn}(d^k_i)\frac{1}{k}\neq1& \mbox{\rm if}\ i<\gamma_x,\vspace{2mm}\\
1+\mbox{\rm sgn}(d^k_i)\frac{1}{k}\neq1& \mbox{\rm if}\ i=\gamma_x,
\end{array}
\right.
\end{equation*}
for $k>K$, where $\mbox{\rm sgn}(\cdot)$ denotes the symbolic function. This contradicts the result that $x^k+\lambda_kd^k\in E$. Thus, $d_i=0$ for any $i\le\gamma_x$. That is $d\in Q$. Thus we have $T^C_E(z)\subseteq Q$. Secondly, we show $T^C_E(z)\supseteq Q$. For any $d\in Q$, define $d^k$ as (\ref{dk-C}). Then for any $z^k\to z$ and $\left\{\lambda_k\right\}\subset E$ with $\lambda_k\to0$, one has $z^k_{\gamma_x}+\lambda_kd^k_{\gamma_x}=1+0=1$ and $z^k_j+\lambda_kd^k_j=0$ for any $j<\gamma_x$ and for all $k\in\mathbb{N}$. Thus, $z^k+\lambda_kd^k\in E$ for all $k\in\mathbb{N}$. So $Q\subseteq T^C_E(z)$. Therefore, $T^C_E(z)=Q$.

Based on the proofs of   (i) and (ii) above, we conclude  that (\ref{NEz}) holds.
\ep

The following definition of M-stationary point of (\ref{originalmodel})  is in line with Definition 4.6 in \cite{BKS2016}, which is similar to the one in \cite{Ye2005} for MPEC (mathematical programs with equilibrium constraints).
\begin{definition}\label{def-M}
$z^*$ is an M-stationary point of (\ref{originalmodel}) if there exists $w^*\in\mathbb{R}^{m+n}$ satisfying $(z^*,w^*)\in\bar{F}$ and $\mu\in\mathbb{R}^{m+n}$ such that the two conditions below hold:
\begin{itemize}
\item $-\nabla f(z^*)-\mu\in N^\#_{E}(z^*)$.
\item $\mu_i=0$, while $z^*_i\neq0$.
\end{itemize}
\end{definition}

\begin{theorem}\label{thm-M-1}
Let $z^*=(x^*,y^*)\in F,$ where $F$ is defined by (\ref{sprbl-set}), and  let $\Gamma$ be defined by (\ref{supphuang}). Then $z^*$ is an M-stationary point of (\ref{originalmodel}) if and only if $\nabla f(z^*)\in\mathbb{R}^{m+n}$ satisfies
\begin{equation}\label{M-stp}
\nabla_if(z^*)=0,\ \forall i\in\Gamma^*\setminus\left\{\gamma_{x^*}\right\}.
\end{equation}
\end{theorem}

\noindent{\bf Proof}. Choose $w^*:=(u^*,v^*)\in\mathbb{R}^{m+n}$ such that
\begin{equation*}
w^*\left\{
\begin{array}{ll}
=0& \mbox{\rm if}\ i\in\Gamma^*,\\
\in(0,1]& {\rm otherwise}.
\end{array}
\right.
\end{equation*}
Then  ${z^*}^\top w^*=0$. Clearly,  $(z^*,w^*)\in\bar{F}$ due to $z^*\in F$. Thus, by Definition \ref{def-M}, $z^*$ is an M-stationary point if and only if there exists $\mu\in\mathbb{R}^{m+n}$ such that $\mu_i=0$ for any $i\in\Gamma^*$ and
 $$
-\nabla f(z^*)-\mu\in N_E(z^*).$$
 By Lemma \ref{N-E-z}, this implies that $z^*$ is an M-stationary point if and only if (\ref{M-stp}) holds.
\ep

The following corollary indicates the relation of M-stationarity and ${\rm N^C}$-stationarity of (\ref{originalmodel}) and also claims that the M-stationarity is  a necessary first-order optimality condition for  (\ref{originalmodel}) .

\begin{corollary}\label{M-N-C}
(i) $z^*=(x^*,y^*)$ is an M-stationary point of (\ref{originalmodel}) if and only if $z^*$ is an ${\rm N^C}$-stationary point of (\ref{originalmodel}).  (ii) $z^*=(x^*,y^*)$ is an optimal solution of (\ref{originalmodel}), then $z^*$ is an M-stationary point.
\end{corollary}

 Item (i) above follows directly from Theorems \ref{thm-N-C} and \ref{thm-M-1}, and Item (ii) follows directly from Theorems \ref{M-N-C} and \ref{B-to-C}. Finally, we summarize the relationship between the stationary points discussed in this work as follows.

\begin{theorem}
The  stationary points of (\ref{originalmodel})  satisfy the following relationship:
\begin{equation*}
\begin{array}{ccccccc}
&&{\rm CW minimum}&&&&\\
&&\Downarrow&&&&\\
L\text{-stationarity}&\Longrightarrow&{\rm N^B}\text{-stationarity}&\Longrightarrow &{\rm N^C}\text{-stationarity}&\Longleftrightarrow &\text{M-stationarity}\\
 & &\Updownarrow& & \Updownarrow& &\\
 & & {\rm T^B}\text{-stationarity}&\Longrightarrow & {\rm T^C}\text{-stationarity}& &
\end{array}
\end{equation*}
\end{theorem}

\section{Conclusion}\label{sect4}

In this paper,  we investigated the first-order optimality conditions for the sparse bilinear least squares problems.  Via some variational geometric tools, we characterized T-type and N-type stationary points of the problem in terms of  tangent  and normal cones in Bouligand sense and Clarke sense. We also showed that an optimal solution to a sparse bilinear least squares problem is necessarily a T-type and/or N-type stationary point in  Bouligand sense and Clarke sense. In addition, the definition of coordinate-wise minima is modified and characterized according to the feasible set, and it is proved that any optimal solution of the considered problem must be a CW minimum. Moreover, we introduced the concept of L-like stationary point based on the so-called like-projection. Based on this concept, we established the first-order necessary optimality condition for the sparse bilinear least squares problem. By reformulating the  problem as a complementarity-type model, we introduced and characterized the M-stationary point, and showed that M-stationarity is a necessary optimality condition for the problem. Furthermore,  the relationship between the various stationary points is also discussed in this work. We believe that the theoretical results established here provide a  fundamental  basis for the development of certain numerical methods for the sparse bilinear least squares problems, which would be an interesting and worthwhile research topic in the near future.

\end{document}